\documentclass[12pt]{article}
\usepackage{amsmath,amssymb,amsthm,latexsym,amscd}

\title{Algebras of distributions \\
for binary semi-isolating formulas \\ of a complete
theory\footnote{\indent \ \ \ The work is
supported by RFBR (grant 12-01-00460-a).}}
\author{Sergey V.
Sudoplatov\footnote{sudoplat@math.nsc.ru}}
\date{October 15, 2012}
\begin{document}
\maketitle

\begin{abstract}
We define a class of algebras describing links of binary
semi-isolating formulas on a set of realizations for a family of
$1$-types of a complete theory. These algebras include algebras of
isolating formulas considered before. We prove that a set of
labels for binary semi-isolating formulas on a set of realizations
for a $1$-type $p$ forms a monoid of a special form with a partial
order inducing ranks for labels, with set-theoretic operations,
and with a composition. We describe the class of these structures.
A description of the class of structures relative to families of
$1$-types is given.

\end{abstract}

In \cite{SuLP}, a series of constructions is introduced admitting
to realize key properties of countable theories and to obtain a
classification of countable models of small (in particular, of
Ehrenfeucht) theories with respect to two basic characteristics:
Rudin--Keisler preorders and distribution functions for numbers of
limit models. The construction of these theories is essentially
based on the definition of special directed graphs with colored
vertices and arcs as well as on the definition of $(n+1)$-ary
predicates that turn prime models over realizations of $n$-types
to prime models over realizations of $1$-types and reducing links
between prime models over finite sets to links between prime
models over elements such that these links are given by definable
sets of arcs and edges.

In the paper, we develop a general approach to the description of
binary links between realizations of $1$-types in terms of labels
of pairwise non-equivalent isolating formulas \cite{ShS} to sets
of labels of semi-isolating formulas.

We use the standard relation algebraic, model-theoretical,
semigroup, and graph-theoretic terminology \cite{HiHo}--\cite{SO2}
as well as some notions, notations, and constructions in
\cite{SuLP, ShS}.

The author thanks Evgeniy A.~Palyutin and Vadim G.~Puzarenko for
useful remarks.

\medskip
\centerline{\bf 1. Preliminary notions, notations, and properties}

\medskip
{\bf Definition} \cite{SuLP, ShS, Su08, BSV}. Let $T$ be a
complete theory, $\mathcal{M}\models T$. Consider types
$p(x),q(y)\in S(\varnothing)$, realized in $\mathcal{M}$, and all
{\em $(p,q)$-preserving}\index{Formula!$(p,q)$-preserving} {\em
$(p,q)$-semi-isolating}\index{Formula!$(p,q)$-semi-isolating},
{\em $(p\rightarrow q)$-}\index{$(p\rightarrow q)$-formulas}, or
{\em $(q\leftarrow p)$-formulas}\index{$(q\leftarrow p)$-formula}
$\varphi(x,y)$ of $T$, i.~e., formulas for which there is $a\in M$
such that $\models p(a)$ and $\varphi(a,y)\vdash q(y)$. Now, for
each such a formula $\varphi(x,y)$, we define a binary relation
$R_{p,\varphi,q}\rightleftharpoons\{(a,b)\mid\mathcal{M}\models
p(a)\wedge\varphi(a,b)\}.$ If $(a,b)\in R_{p,\varphi,q}$, $(a,b)$
is called a {\em $(p,\varphi,q)$-arc}\index{$(p,\varphi,q)$-arc}.
If $\varphi(a,y)$ is principal (over $a$), the $(p,\varphi,q)$-arc
$(a,b)$ is also {\em
principal}\index{$(p,\varphi,q)$-arc!principal}\index{Arc!principal}.

If,~in~addition, $\varphi(x,y)$ is a {\em $(p\leftrightarrow
q)$-formula}\index{$(p\leftrightarrow q)$-formula}, i.~e., it is
both a $(p\rightarrow q)$- and a $(q\rightarrow p)$-formula then
the set $[a,b]\rightleftharpoons\{(a,b),(b,a)\}$ is said to be a
{\em $(p,\varphi,q)$-edge}\index{$(p,\varphi,q)$-edge}. If the
$(p,\varphi,q)$-edge $[a,b]$ consists of principal
$(p,\varphi,q)$- and $(q,\varphi(y,x),p)$-arcs then $[a,b]$ is a
{\em principal}
$(p,\varphi,q)$-edge\index{$(p,\varphi,q)$-edge!principal}.

$(p,\varphi,q)$-arcs and $(p,\varphi,q)$-edges are called {\em
arcs}\index{Arc} and {\em edges}\index{Edge} respectively if we
say about fixed or some $(p\rightarrow q)$-formula $\varphi(x,y)$.
If $(a,b)$ is a principal $(p,\varphi,q)$-arc such that the pair
$(b,a)$ is not a principal arc (on any formula), that is
$(b,a)\notin{\rm SI}_{\{p,q\}}$, then $(a,b)$ is called {\em
irreversible}\index{Arc!principal!irreversible}. If $(a,b)$ is a
$(p,\varphi,q)$-arc and $(b,a)$ is not a $(q,\varphi,p)$-arc then
$(a,b)$ is also an {\em irreversible} arc.\index{Arc!irreversible}

\medskip
For types $p(x),q(y)\in S(\varnothing)$, we denote by ${\rm
SICF}(p,q)$\index{${\rm SICF}(p,q)$} the set of $(p\rightarrow
q)$-formulas $\varphi(x,y)$ such that $\{\varphi(a,y)\}$ is
consistent for $\models p(a)$. Let ${\rm SICE}(p,q)$\index{${\rm
SICE}(p,q)$} be the set of pairs of formulas
$(\varphi(x,y),\psi(x,y))\in{\rm SICF}(p,q)$ such that for any
(some) realization $a$ of $p$ the sets of solutions for
$\varphi(a,y)$ and $\psi(a,y)$ coincide. Clearly, ${\rm
SICE}(p,q)$ is an equivalence relation on the set ${\rm
SICF}(p,q)$. Notice that each ${\rm SICE}(p,q)$-class $E$
corresponds to either a set of $(p,\varphi,q)$-edges, or a set of
irreversible $(p,\varphi,q)$-arcs, or simultaneously a set of
$(p,\varphi,q)$-edges and of irreversible $(p,\varphi,q)$-arcs
linking realizations of $p$ and $q$ by any (some) formula
$\varphi$ in $E$.  Thus the quotient ${\rm SICF}(p,q)/{\rm
SICE}(p,q)$ is represented as a disjoint union of sets ${\rm
SICFE}(p,q)$\index{${\rm SICFE}(p,q)$}, ${\rm
SICFA}(p,q)$\index{${\rm SICFA}(p,q)$}, and ${\rm
SICFM}(p,q)$\index{${\rm SICFM}(p,q)$}, where ${\rm SICFE}(p,q)$
consists of ${\rm SICE}(p,q)$-classes correspondent to sets of
edges, ${\rm SICFA}(p,q)$ consists of ${\rm SICE}(p,q)$-classes
correspondent to sets of irreversible arcs, and ${\rm SICFM}(p,q)$
consists of ${\rm SICE}(p,q)$-classes correspondent to sets
containing edges and irreversible arcs.

The sets ${\rm SICF}(p,p)$, ${\rm SICE}(p,p)$, ${\rm SICFE}(p,p)$,
${\rm SICFA}(p,p)$, and ${\rm SICFM}(p,p)$ are denoted by ${\rm
SICF}(p)$,\index{${\rm SICF}(p)$} ${\rm SICE}(p)$,\index{${\rm
SICE}(p)$} ${\rm SICFE}(p)$\index{${\rm SICFE}(p)$}, ${\rm
SICFA}(p)$\index{${\rm SICFA}(p)$}, and ${\rm
SICFM}(p)$\index{${\rm SICFM}(p)$} respectively.

Let $T$ be a complete theory without finite models,
$U=U^-\,\dot{\cup}\,\{0\}\,\dot{\cup}\,U^+\,\dot{\cup}\,U'$ be an
alphabet of cardinality $\geq|S(T)|$ and consisting of {\em
negative elements}\index{Element!negative} $u^-\in
U^-$,\index{$U^-$} {\em positive elements}\index{Element!positive}
$u^+\in U^+$\index{$U^+$}, {\em neutral
elements}\index{Element!neutral} $u'\in U'$\index{$U'$}, and zero
$0$. As usual, we write $u<0$ for any $u\in U^-$ and $u>0$ for any
$u\in U^+$. The set $U^-\cup\{0\}$ is denoted by $U^{\leq
0}$\index{$U^{\leq 0}$} and $U^+\cup\{0\}$ is denoted by $U^{\geq
0}$\index{$U^{\geq 0}$}. Elements of $U$ are called {\em
labels}.\index{Label}

Let $\nu(p,q)\mbox{\rm : }{\rm SICF}(p,q)/{\rm SICE}(p,q)\to U$ be
injective {\em labelling functions},\index{Function!labelling}
$p(x),q(y)\in S(\varnothing)$, for which negative elements
correspond to the classes in ${\rm SICFA}(p,q)/{\rm SICE}(p,q)$,
positive elements and $0$ correspond to the classes in ${\rm
SICFE}(p,q)/{\rm SICE}(p,q)$ such that $0$ is defined only for
$p=q$ and is represented by the formula $(x\approx y)$, and
neutral elements code the classes in ${\rm SICFM}(p,q)/{\rm
SICE}(p,q)$, $\nu(p)\rightleftharpoons\nu(p,p)$\index{$\nu(p)$}.
We additionally assume that $\rho_{\nu(p)}\cap\rho_{\nu(q)}=\{0\}$
for $p\ne q$ where, as usual, we denote by $\rho_f$ the image of
the function $f$) and
$\rho_{\nu(p,q)}\cap\rho_{\nu(p',q')}=\varnothing$ if $p\ne q$ and
$(p,q)\ne(p',q')$. Labelling functions with the properties above
as well families of these functions are said to be {\em
regular}\index{Function!labelling!regular}\index{Family!regular}.
Further we shall consider only regular labelling functions and
their regular families.

The labels, correspondent to isolating formulas, are said to be
{\em isolating}\index{Label!isolating} whereas each label in
$\bigcup\limits_{p,q\in S^1(\varnothing)}\rho_{\nu(p,q)}$ is {\em
semi-isolating}\index{Label!semi-isolating}. By the definition,
each isolating label belongs to
$U^-\,\dot{\cup}\,\{0\}\,\dot{\cup}\,U^+$, i.~e., it is not
neutral.

We denote by $\theta_{p,u,q}(x,y)$\index{$\theta_{p,u,q}(x,y)$}
formulas in ${\rm SICF}(p,q)$ with a label $u\in\rho_{\nu(p,q)}$.
If the type $p$ is fixed and $p=q$ then the formula
$\theta_{p,u,q}(x,y)$ is denoted by
$\theta_u(x,y)$.\index{$\theta_u(x,y)$}

Note that if $\theta_{p,u,q}(x,y)$ and $\theta_{q,v,p}(x,y)$ are
formulas witnessing that for realizations $a$ and $b$ of $p$ and
$q$ respectively the pairs $(a,b)$ and $(b,a)$ belong to ${\rm
SI}_{\{p,q\}}$, then the formula
$\theta_{p,u,q}(x,y)\wedge\theta_{q,v,p}(y,x)$ witnesses that
$[a,b]$  is a $(p,\varphi,q)$-edge. If the edge $[a,b]$ is
principal and $\theta_{p,u,q}(a,y)$ is an isolating formula such
that $\models\theta_{p,u,q}(a,b)$, $\models p(a)$, then the label
$u$ is {\em invertible}\index{Label!invertible} and a label $v\in
U^{\geq 0}$ corresponds uniquely to $u$ such that
$\theta_{q,v,p}(b,y)$ is an isolating formula with
$\models\theta_{q,v,p}(b,a)$, and vice versa. The labels $u$ and
$v$ are {\em reciprocally inverse}\index{Labels!reciprocally
inverse} and are denoted by $v^{-1}$\index{$v^{-1}$} and
$u^{-1}$\index{$u^{-1}$} respectively. In general case, each label
$u\in U^{\geq 0}$ has a (nonempty) set of {\em inverse}
labels\index{Label!inverse}\index{Set!of inverse labels} in
$U^{\geq 0}$, denoted also by $u^{-1}$\index{$u^{-1}$}. Note that
independently on a label $u\in U^{\geq 0}$, for which a formula
$\theta_{p,u,q}(x,y)$ witnesses that $[a,b]$ is a
$(p,\varphi,q)$-edge, the set $u^{-1}$ includes all labels $v\in
U^{\geq 0}$ such that $[b,a]$ is a $(q,\theta_{q,v,p},p)$-edge.

Neutral labels correspond, for instance, the formulas
$\theta_{p,u,q}(x,y)\vee\theta_{p,v,q}(x,y)$, where $u<0$ and
$v\geq 0$.

For types $p_1,p_2,\ldots,p_{k+1}\in S^1(\varnothing)$ and sets
$X_1,X_2,\ldots,X_k\subseteq U$ of labels we denote by\index{${\rm
SI}(p_1,X_1,p_2,X_2,\ldots,p_k,X_k,p_{k+1})$}
$${\rm
SI}(p_1,X_1,p_2,X_2, \ldots, p_k,X_k,p_{k+1})$$ the set of all
labels $u\in U$ correspondent to formulas
$\theta_{p_1,u,p_{k+1}}(x,y)$ satisfying, for realizations $a$ of
$p_1$ and some $u_1\in X_1,\ldots,u_k\in X_k$, the following
condition:
$$
\theta_{p_1,u,p_{k+1}}(a,y)\vdash\theta_{p_1,u_1,p_2,u_2,\ldots,p_k,u_k,p_{k+1}}(a,y),$$
where\index{$\theta_{p_1,u_1,p_2,u_2,\ldots,p_k,u_k,p_{k+1}}(x,y)$}
$$\theta_{p_1,u_1,p_2,u_2,\ldots,p_k,u_k,p_{k+1}}(x,y)\rightleftharpoons$$
$$\rightleftharpoons\exists x_2,x_3,\ldots
x_{k-1}(\theta_{p_1,u_1,p_2}(x,x_2)\wedge\theta_{p_2,u_2,p_3}(x_2,x_3)\wedge\ldots$$
$$\ldots\wedge\theta_{p_{k-1},u_{k-1},p_k}(x_{k-1},x_k)\wedge\theta_{p_k,u_k,p_{k+1}}(x_k,y)).
$$
In view of transitivity of semi-isolation, each formula
$\theta_{p_1,u_1,p_2,u_2,\ldots,p_k,u_k,p_{k+1}}(x,y)$ has a label
in $\rho_{\nu(p_1,p_{k+1})}$.

Thus the Boolean $\mathcal{P}(U)$ of $U$ is the universe of an
{\em algebra $\mathfrak{A}$ of distributions of binary
semi-isolating formulas}\index{Algebra!of distributions of binary
semi-isolating formulas} with $k$-ary operations
$${\rm
SI}(p_1,\cdot,p_2,\cdot, \ldots, p_k,\cdot,p_{k+1}),$$ where
$p_1,\ldots,p_{k+1}\in S^1(\varnothing)$. This algebra has a
natural restriction to any family $R\subseteq S^1(\varnothing)$ as
well as to the algebras of distributions of binary {\em isolating}
formulas \cite{ShS}. Besides, if $U_0$ is a subalphabet of $U$
then the restriction of the universe of $\mathfrak{A}$ to the set
$\mathcal{P}(U_0)$ and  the restrictions for values of operations
to the set $U_0$ forms, possibly partial, algebra
$\mathfrak{A}\upharpoonright
U_0$.\index{$\mathfrak{A}\upharpoonright U_0$}

Note that if some set $X_i$ is disjoint with
$\rho_{\nu(p_i,p_{i+1})}$, in particular, if it is empty then
$${\rm
SI}(p_1,X_1,p_2,X_2,\ldots,p_k,X_k,p_{k+1})=\varnothing,$$  and if
each $X_i$ has common elements with $\rho_{\nu(p_i,p_{i+1})}$ then
$${\rm
SI}(p_1,X_1,p_2,X_2,\ldots,p_k,X_k,p_{k+1})\ne\varnothing.$$

Note also that if $X_i\not\subseteq \rho_{\nu(p_i,p_{i+1})}$ for
some $i$ then $${\rm SI}(p_1,X_1,p_2,X_2, \ldots,
p_k,X_k,p_{k+1})=$$
$$= {\rm
SI}(p_1,X_1\cap\rho_{\nu(p_1,p_2)},p_2,X_2\cap\rho_{\nu(p_2,p_3)},
\ldots, p_k,X_k\cap\rho_{\nu(p_k,p_{k+1})},p_{k+1}).$$

In view of the previous equation, further, considering values
$${\rm
SI}(p_1,X_1,p_2,X_2,\ldots,p_k,X_k,p_{k+1}),$$ we shall assume
that $X_i\subseteq \rho_{\nu(p_i,p_{i+1})}$, $i=1,\ldots,k$.

If each set $X_i$ is a singleton consisting of an element $u_i$
then we use $u_i$ instead of $X_i$ in ${\rm
SI}(p_1,X_1,p_2,X_2,\ldots,p_k,X_k,p_{k+1})$ and write\index{${\rm
SI}(p_1,u_1,p_2,u_2,\ldots,p_k,u_k,p_{k+1})$}
$${\rm
SI}(p_1,u_1,p_2,u_2,\ldots,p_k,u_k,p_{k+1}).$$

By the definition the following equality holds:
$${\rm
SI}(p_1,X_1,p_2,X_2, \ldots, p_k,X_k,p_{k+1})=$$
$$=\cup\{{\rm
SI}(p_1,u_1,p_2,u_2, \ldots, p_k,u_k,p_{k+1})\mid u_1\in
X_1,\ldots, u_k\in X_k\}.$$ Hence the specification of ${\rm
SI}(p_1,X_1,p_2,X_2, \ldots, p_k,X_k,p_{k+1})$ is reduced to the
specifications of ${\rm SI}(p_1,u_1,p_2,u_2, \ldots,
p_k,u_k,p_{k+1})$. Note also that ${\rm SI}(p,X,q)=X$ for any
$X\subseteq\rho_{\nu(p,q)}$.

Clearly, if $u_i=0$ then $p_i=p_{i+1}$ for nonempty sets
$${\rm
SI}(p_1,u_1,p_2,u_2,\ldots,p_i,0,p_{i+1},\ldots,
p_k,u_k,p_{k+1})$$ and the following conditions hold:
$$
{\rm SI}(p_1,0,p_1)=\{0\},
$$
$${\rm
SI}(p_1,u_1,p_2,u_2,\ldots,p_i,0,p_{i+1},\ldots,
p_k,u_k,p_{k+1})=$$
$$={\rm
SI}(p_1,u_1,p_2,u_2,\ldots,p_i,u_{i+1},p_{i+2},\ldots,
p_k,u_k,p_{k+1}).$$

If all types $p_i$ equal to a type $p$ then we write ${\rm
SI}_p(X_1,X_2,\ldots,X_k)$\index{${\rm SI}_p(X_1,X_2,\ldots,X_k)$}
and ${\rm SI}_p(u_1,u_2,\ldots,u_k)$\index{${\rm
SI}_p(u_1,u_2,\ldots,u_k)$} as well as $\lceil
X_1,X_2,\ldots,X_k\rceil_p$\index{$\lceil X_1,X_2,\ldots,X_k\rceil
_p$} and $\lceil u_1,u_2,\ldots,u_k\rceil _p$\index{$\lceil
u_1,u_2,\ldots,u_k\rceil _p$} instead of
$${\rm
SI}(p_1,X_1,p_2,X_2,\ldots,p_k,X_k,p_{k+1})$$ and
$${\rm
SI}(p_1,u_1,p_2,u_2,\ldots,p_k,u_k,p_{k+1})$$ respectively. We
omit the index ${\cdot}_p$ if the type $p$ is fixed. In this case,
we write
$\theta_{u_1,u_2,\ldots,u_k}(x,y)$\index{$\theta_{u_1,u_2,\ldots,u_k}(x,y)$}
instead of $\theta_{p,u_1,p,u_2,\ldots,p,u_k,p}(x,y)$.

\medskip
{\bf Proposition 1.1.} {\em $(1)$ If $p,q\in S^1(T)$ are principal
types then $\rho_{\nu(p,q)}\cup\rho_{\nu(q,p)}\subseteq U^{\geq
0}$.

$(2)$ If $p,q\in S^1(T)$, $p$ is a principal type and $q$ is a
non-principal type then $\rho_{\nu(p,q)}=\varnothing$ and
$\rho_{\nu(q,p)}\subseteq U^-$.}

\medskip
{\em Proof.} (1) If $\rho_{\nu(p,q)}$ contains a label $u\notin
U^{\geq 0}$ then there are realizations $a$ and $b$ of $p$ and $q$
respectively such that $(a,b)\in{\rm SI}_{\{p,q\}}$ and
$(b,a)\notin{\rm SI}_{\{p,q\}}$. But since $p(x)$ contains a
principal formula $\varphi(x)$, this formula witnesses that
$(b,a)\in{\rm SI}_{\{p,q\}}$. The contradiction implies that
$\rho_{\nu(p,q)}\subseteq U^{\geq 0}$. Similarly we obtain
$\rho_{\nu(q,p)}\subseteq U^{\geq 0}$.

(2) Let $\varphi(x)$ be a principal formula of $p(x)$. If $\models
p(a)$, $\models q(b)$, and $(a,b)\in{\rm SI}_{\{p,q\}}$ that
witnessed by a formula $\theta_u(x,y)$, the formula $\exists
x(\varphi(x)\wedge\theta_u(x,y))$ isolates $q(y)$. Since $q$ is
not isolated we obtain $\rho_{\nu(p,q)}=\varnothing$. By the same
reason, $\rho_{\nu(q,p)}\subseteq U^-$.~$\Box$

\medskip
{\bf Corollary 1.2.} {\em If $p(x)$ is a principal type then
$\rho_{\nu(p)}\subseteq U^{\geq 0}$.}

\medskip
{\bf Proposition 1.3.} {\em Let $p_1,p_2,\ldots,p_{k+1}$ be types
in $S^1(\varnothing)$. The following assertions hold.

$(1)$ If $u_i\in\rho_{\nu(p_i,p_{i+1})}$, $i=1,\ldots,k$, and some
$u_i$ is negative then
$${\rm
SI}(p_1,u_1,p_2,u_2, \ldots, p_k,u_k,p_{k+1})\subseteq U^-.$$

$(2)$ If $u_i\in\rho_{\nu(p_i,p_{i+1})}\cap U^{\geq 0}$,
$i=1,\ldots,k$, then
$${\rm
SI}(p_1,u_1,p_2,u_2,\ldots,p_k,u_k,p_{k+1})\subseteq U^{\geq 0}.$$

$(3)$ If $u_i\in\rho_{\nu(p_i,p_{i+1})}\cap (U^{\geq 0}\cup U')$,
$i=1,\ldots,k$, and some $u_i$ belongs to  $U'$ then
$${\rm
SI}(p_1,u_1,p_2,u_2,\ldots,p_k,u_k,p_{k+1})\subseteq U'.$$

$(4)$ If $u_i\in\rho_{\nu(p_i,p_{i+1})}\cap U^{\geq 0}$,
$i=1,\ldots,k$, then all elements of the set $X\rightleftharpoons
{\rm SI}(p_1,u_1,p_2,u_2,\ldots,p_k,u_k,p_{k+1})$ are invertible
and the set $X^{-1}\rightleftharpoons\cup\{v^{-1}\mid v\in
X\}$\index{$X^{-1}$} coincides with the set ${\rm
SI}(p_{k+1},u^{-1}_k,p_k,u^{-1}_{k-1},\ldots,p_2,u^{-1}_1,p_{1})$.}

\medskip
{\em Proof.} (1)--(3) follow by the transitivity of
semi-isolation.

(4) All elements in $X$ are invertible by (2). Let $v'$ be an
element in $v^{-1}\subseteq X^{-1}$. Then for any
$(p_1,\theta_{p_1,v,p_{k+1}},p_{k+1})$-edge $[a,b]$ the following
conditions hold:

(a) there are realizations $a_i$ of $p_i$, $i=1,\ldots,k+1$, such
that $a_0=a$, $a_{k+1}=b$,
$\models\theta_{p_i,u_i,p_{i+1}}(a_i,a_{i+1})$, $i=1,\ldots,k$;

(b) $[b,a]$ is a $(p_{k+1},\theta_{p_{k+1},v',p_1},p_1)$-edge.

Since $[a_{i+1},a_{i}]$ is an $u'_i$-edge for any $u'_i\in
u_i^{-1}$, $i=1,\ldots,k$, then
$$
\theta_{p_{k+1},v',p_1}(b,x)\vdash\theta_{p_{k+1},u'_k,p_k,u'_{k_1}\ldots,p_2,u'_1,p_1}(b,x),
$$
whence, $v'\in{\rm
SI}(p_{k+1},u^{-1}_k,p_k,u^{-1}_{k-1},\ldots,p_2,u^{-1}_1,p_{1})$.

If $v'\in{\rm
SI}(p_{k+1},u'_k,p_k,u'_{k-1},\ldots,p_2,u^{-1}_1,p_{1})$,
$u'_i\in u^{-1}_i$, $i=1,\ldots,k$, then $v'\geq 0$ and for any
$(p_{k+1},\theta_{p_{k+1},v',p_1},p_1)$-edge $[b,a]$ there are
realizations $b_i$ of $p_i$, $i=1,\ldots,k+1$, such that
$b_{k+1}=b$, $b_{1}=a$,
$\models\theta_{p_{i+1},u'_i,p_{i}}(b_{i+1},b_{i})$,
$i=1,\ldots,k$. We have
$\models\theta_{p_i,u_i,p_{i+1}}(b_i,b_{i+1})$, $i=1,\ldots,k$,
and so the elements $b_1,\ldots,b_{k+1}$ witness that $[a,b]$ is a
$(p_1,\theta_{p_{1},u_1,p_2,u_{2}\ldots,p_k,u_k,p_{k+1}},p_{k+1})$-edge.
Then $v'$ belongs to $v^{-1}$, where $v\in X$ is a label for the
formula $\theta_{p_{1},u_1,p_2,u_{2}\ldots,p_k,u_k,p_{k+1}}$.
Thus, $v'\in X^{-1}$.~$\Box$

\medskip
{\bf Corollary 1.4.} {\em Restrictions of  $U$ to the sets
$U^{\leq 0}$, $U^{\geq 0}$, and $U^{\geq 0}\cup U'$ form
subalgebras of the algebra of distributions of binary
semi-isolating formulas. The operation of inversion is coordinated
with the operations of the algebra.}

\medskip
\centerline{\bf 2. Preordered algebras} \centerline{\bf of
distributions of binary semi-isolating formulas}
\medskip

For the set $U$ of labels in the algebra $\mathfrak{A}$ of binary
semi-isolating formulas of theory $T$, we define the following
relation $\unlhd$: if $u,v\in U$ then $u\unlhd v$\index{$u\unlhd
v$} if and only if $u=v$, or $u,v\in\rho_{\nu(p,q)}$ for some
types $p,q\in S^1(\varnothing)$ and
$\theta_{p,u,q}(a,y)\vdash\theta_{p,v,q}(a,y)$ for some (any)
realization $a$ of $p$. If $u\unlhd v$ and $u\ne v$ we write
$u\lhd v$.\index{$u\lhd v$}

By the definition the relation $\unlhd$ is reflexive and
transitive. It is antisymmetric since distinct labels correspond
to non-equivalent formulas.

Below we consider some properties for the substructures $\langle
\rho_{\nu(p,q)};\unlhd\rangle$ of the partially ordered set
$\langle U;\unlhd\rangle$.

\medskip
{\bf Proposition 2.1.} {\em $(1)$ For any types $p,q\in
S^1(\varnothing)$ the partially ordered set $\langle
\rho_{\nu(p,q)};\unlhd\rangle$ forms a upper semilattice.

$(2)$ An element $u\in \rho_{\nu(p,q)}$ is $\unlhd$-minimal if and
only if for a realization $a$ of $p$, the formula
$\theta_{p,u,q}(a,y)$ is isolating.

{\rm $(3)$ (monotony).}\index{Monotony} If $u,v\in
\rho_{\nu(p,q)}$ and $u\unlhd v$ then $v\in U^\delta$,
$\delta\in\{-,+\}$, implies $u\in U^\delta$, and if $u\in U'$ then
$v\in U'$.}

\medskip
{\em Proof.} (1) If $u_1,u_2\in\rho_{\nu(p,q)}$ then for the
formulas $\theta_{p,u_1,q}(x,y)$ and $\theta_{p,u_2,q}(x,y)$ the
label $v$ for the formula
$\theta_{p,u_1,q}(x,y)\vee\theta_{p,u_2,q}(x,y)$ is a supremum for
the labels $u_1$ and $u_2$.

(2) If $\theta_{p,u,q}(a,y)$ is an isolating formula then the
label $u$ is $\unlhd$-minimal by the definition. If the formula
$\theta_{p,u,q}(a,y)$ is not isolating then there is a formula
$\varphi(a,y)$ such that the semi-isolating formulas
$\theta_{p,u,q}(a,y)\wedge\varphi(a,y)$ and
$\theta_{p,u,q}(a,y)\wedge\neg\varphi(a,y)$ are consistent. For
the labels $v_1$ and $v_2$ of these formulas, we have $v_1\ne
v_2$, $v_1\lhd u$, and $v_2\lhd u$.

(3) If $v\in\rho_{\nu(p,q)}\cap U^-$ then for any solution  $b$ of
the formula $\theta_{p,v,q}(a,y)$, where $\models p(a)$, the pair
$(a,b)$ is an irreversible arc. Hence, for any solution $b$ of
$\theta_{p,u,q}(a,y)$, where $u\unlhd v$, the pair $(a,b)$ is also
an irreversible arc and so $u$ belongs to $U^-$. Replacing arcs by
edges, the same arguments show that $u\unlhd v$ and $v\in U^+$
imply $u\in U^+$. If $u\in U'$ then the set of pairs $(a,b)$ for
the formula $\theta_{p,u,q}(a,y)$ contains both irreversible and
reversible arcs. This property is preserved for any label $v$ with
$u\unlhd v$, whence $v\in U'$.~$\Box$

\medskip
The partial order $\unlhd$ has a natural extension to a preorder
on the set $\mathcal{P}(U)$: for any sets $X,Y\in\mathcal{P}(U)$
we put $X\unlhd Y$\index{$X\unlhd Y$} if $X=\varnothing$, or for
any $x\in X$ there is $y\in Y$ with $x\unlhd y$ and for any $y\in
Y$ there is $x\in X$ with $x\unlhd y$. Thus, the algebra
$\mathfrak{A}$ is transformed to the preordered algebra
$\langle\mathfrak{A};\unlhd\rangle$ with the monotonic property
with respect its restrictions to the sets $U^{\leq 0}$, $U^{\geq
0}$, and $U'$.

Another natural expansion of the now preordered algebra
$\langle\mathfrak{A};\unlhd\rangle$ is based on the the properties
mentioned that if $u_1,u_2\in\rho_{\nu(p,q)}$ and
$v\in\rho_{\nu(q,r)}$ then the formulas
$\theta_{p,u_1,q,v,r}(x,y)$ and
$\theta_{p,u_1,q}(x,y)\vee\theta_{p,u_2,q}(x,y)$ as well as
$\theta_{p,u_1,q}(x,y)\wedge\theta_{p,u_2,q}(x,y)$ and
$\theta_{p,u_1,q}(x,y)\wedge\neg\theta_{p,u_2,q}(x,y)$ (if the
formulas $\theta_{p,u_1,q}(a,y)\wedge\theta_{p,u_2,q}(a,y)$ and
$\theta_{p,u_1,q}(a,y)\wedge\neg\theta_{p,u_2,q}(a,y)$ are
consistent for $\models p(a)$) have labels in $U$. We denote these
labels by $u_1\circ v$,\index{$u_1\circ v$} $u_1\vee
u_2$,\index{$u_1\vee u_2$} $u_1\wedge u_2$,\index{$u_1\wedge u_2$}
and $u_1\wedge\neg u_2$\index{$u_1\wedge\neg u_2$} respectively.
The last label is also denoted by $\neg u_2\wedge u_1$\index{$\neg
u_2\wedge u_1$}. The label $u_1\circ v$ is the {\em
composition}\index{Composition!of labels} of labels $u_1$ and $v$;
$u_1\vee u_2$ is the {\em union}\index{Union!of labels} or the
{\em disjunction}\index{Disjunction!of labels} of labels $u_1$ and
$u_2$; $u_1\wedge u_2$ is their {\em
intersection}\index{Intersection!of labels} or {\em
conjunction}\index{Conjunction!of labels}; $u_1\wedge\neg u_2$ is
the {\em relative complement}\index{Complement relative} of $u_2$
in $u_1$.

Clearly, $u_1\unlhd u_1\vee u_2$, $u_2\unlhd u_1\vee u_2$,
$u_1\wedge u_2\unlhd u_1$, $u_1\wedge u_2\unlhd u_2$, $u_1\wedge
\neg u_2\unlhd u_1$.

We set
$$(p,(u_1\circ v),r)\rightleftharpoons\left\{\begin{array}{ll}
\{u_1\circ v\},&\mbox{ if }u_1\in\rho_{\nu(p,q)}\mbox{ and
}v\in\rho_{\nu(q,r)},\\
\varnothing,&\mbox{ if }u_1\notin\rho_{\nu(p,q)}\mbox{ or
}v\notin\rho_{\nu(q,r)},\end{array}\right.$$
$$(p,(u_1\vee u_2),q)\rightleftharpoons\left\{\begin{array}{ll}
\{u_1\vee u_2\},&\mbox{ if }u_1\in\rho_{\nu(p,q)}\mbox{ and
}u_2\in\rho_{\nu(p,q)},\\
\{u_1\},&\mbox{ if }u_1\in\rho_{\nu(p,q)}\mbox{ and
}u_2\notin\rho_{\nu(p,q)},\\
\{u_2\},&\mbox{ if }u_1\notin\rho_{\nu(p,q)}\mbox{ and
}u_2\in\rho_{\nu(p,q)},\\ \varnothing,&\mbox{ if
}u_1\notin\rho_{\nu(p,q)}\mbox{ and
}u_2\notin\rho_{\nu(p,q)},\end{array}\right.$$
$$(p,(u_1\wedge u_2),q)\rightleftharpoons\left\{\begin{array}{ll}
\{u_1\wedge u_2\},&\mbox{ if
}u_1\in\rho_{\nu(p,q)},u_2\in\rho_{\nu(p,q)}\\
\,&\mbox{ and
}\models\exists y(\theta_{p,u_1,q}(a,y)\wedge\theta_{p,u_2,q}(a,y)),\\
\varnothing,&\mbox{ otherwise},\end{array}\right.$$
$$(p,(u_1\wedge \neg u_2),q)\rightleftharpoons\left\{\begin{array}{ll}
\{u_1\wedge\neg u_2\},&\mbox{ if
}u_1\in\rho_{\nu(p,q)},u_2\in\rho_{\nu(p,q)}\\
\,&\mbox{ and
}\models\exists y(\theta_{p,u_1,q}(a,y)\wedge\neg\theta_{p,u_2,q}(a,y)),\\
\varnothing,&\mbox{ otherwise},\end{array}\right.$$
$$(p,(X_1\,\tau\, X_2),q)\rightleftharpoons\cup\{(p,(u_1\,\tau\, u_2),q)\mid u_1\in X_1,u_2\in X_2\},\,\,\tau\in\{\circ,\vee,\wedge\},$$
$$(p,(X_1\wedge\neg X_2),q)\rightleftharpoons(p,(\neg X_2\wedge X_1),q)\rightleftharpoons$$
$$\rightleftharpoons\cup\{(p,(u_1\wedge\neg u_2),q)\mid u_1\in X_1,u_2\in X_2\},\,\,X_1,X_2\in\mathcal{P}(U).$$

Labels $u_1$ and $u_2$ are {\em
consistent}\index{Labels!consistent} if $u_1\wedge u_2\in U$. If
$u_1\wedge u_2=\varnothing$ the labels $u_1$ and $u_2$ are called
{\em inconsistent}.\index{Labels!inconsistent}

The preordered algebra $\langle\mathfrak{A};\unlhd\rangle$
equipped with binary operations $(p,(\cdot\,\tau\, \cdot),q)$,
$\tau\in\{\vee,\wedge,\circ\}$, and $(p,(\cdot\,\wedge\,\neg\,
\cdot),q)$, $p,q\in S^1(\varnothing)$, is called a {\em preordered
algebra with relative set-theoretic operations and the
composition}\index{Algebra!preordered with relative set-theoretic
operations and the composition} or briefly a {\em {\rm
POSTC}-algebra}.\index{{\rm POSTC}-algebra}

For any types $p,q\in S^1(\varnothing)$ the structure
$\langle\rho_{\nu(p,q)}\cup\{\varnothing\};\vee,\wedge,\varnothing\rangle$
with operations $\vee$ and $\wedge$ on labels, being extended by
equalities $u\vee\varnothing=u$, $u\wedge\varnothing=\varnothing$,
where $u\in\rho_{\nu(p,q)}\cup\{\varnothing\}$, is an {Ershov
algebra},\index{Algebra!Ershov} i.~e., a {\em distributive lattice
with zero $\varnothing$ and relative
complements}\index{Lattice!distributive with relative complements}
\cite{Er1} such that for any $u,v\in\rho_{\nu(p,q)}$ if $u\unlhd
v$ and $u'=\neg u\wedge v$ is a label then $u\wedge
u'=\varnothing$ and $u\vee u'=v$, and if the label $u'$ does not
exist then $u=v$.

A label $u\in U$ is an {\em atom}\index{Atom} or an {\em atomic
label}\index{Label!atomic} if $u$ is a $\unlhd$-minimal element in
$U$, i.~e., for any label $v\in U$ if $v\unlhd u$ then $v=u$.

By Proposition 2.1, the set of atoms equals the set of isolating
labels and, thus, each atom $u\in\rho_{\nu(p,q)}$ is represented
by an isolated formula $\theta_{p,u,q}(a,y)$, where $\models
p(a)$.

Let $R$ be a nonempty family of types in $S^1(\varnothing)$,
$\mathfrak{A}_R$\index{$\mathfrak{A}_R$} be a restriction of {\rm
POSTC}-algebra $\mathfrak{A}$ to the family $R$. The structure
$\mathfrak{A}_R$ is {\em atomic}\index{Structure!atomic} if for
any types $p,q\in R$ and for any label $u\in\rho_{\nu(p,q)}$ there
is an atom $v\in\rho_{\nu(p,q)}$ such that $v\unlhd u$. The {\rm
POSTC}-algebra $\mathfrak{A}$ is called {\em
$R$-atomic}\index{{\rm POSTC}-algebra!$R$-atomic} if
$\mathfrak{A}_R$ is atomic. If $R=S^1(\varnothing)$ then the
$R$-atomic {\rm POSTC}-algebra is called {\em atomic}.\index{{\rm
POSTC}-algebra!atomic}

Using the definition of atomic structure, of $R$-atomic {\rm
POSTC}-algebra, and of small theory we obtain the following
assertions.

\medskip
{\bf Proposition 2.2.} {\em If $R$ is a nonempty family of types
in $S^1(\varnothing)$ and for any type $p\in R$, there is an
atomic model $\mathcal{M}_p$ over a realization of $p$, then the
{\rm POSTC}-algebra $\mathfrak{A}$ is $R$-atomic.}

\medskip
{\bf Corollary 2.3.} {\em If $T$ is a small theory then the {\rm
POSTC}-algebra $\mathfrak{A}$ is atomic.}

\medskip
\centerline{\bf 3. Ranks and degrees of semi-isolation}
\medskip

The following definition is a local variation of Morley rank
\cite{Mo65}.

\medskip
{\bf Definition.} For triples $(p,u,q)$, where $p,q\in
S^1(\varnothing)$, $u\in U\cup\{\varnothing\}$, we define
inductively the {\em rank ${\rm si}(p,u,q)$ of
semi-isolation}\index{Rank!of semi-isolation}:

(1) ${\rm si}(p,u,q)=0$ if $u\notin\rho_{\nu(p,q)}$;

(2) ${\rm si}(p,u,q)\geq 1$ if $u\in\rho_{\nu(p,q)}$;

(3) for a positive ordinal $\alpha$, ${\rm si}(p,u,q)\geq\alpha+1$
if there is a set $\{v_i\mid i\in\omega\}$ of pairwise
inconsistent labels such that $v_i\lhd u$ and ${\rm
si}(p,v_i,q)\geq\alpha$, $i\in\omega$;

(4) for a limit ordinal $\alpha$, ${\rm si}(p,u,q)\geq\alpha$ if
${\rm si}(p,u,q)\geq\beta$ for any $\beta\in\alpha$.

As usual, we write ${\rm si}(p,u,q)=\alpha$ if ${\rm
si}(p,u,q)\geq\alpha$ and ${\rm si}(p,u,q)\not\geq\beta$ for
$\alpha\in\beta$; ${\rm si}(p,u,q)\rightleftharpoons\infty$ if
${\rm si}(p,u,q)\geq\alpha$ for any ordinal $\alpha$.

If types $p$ and $q$ are fixed, we write ${\rm si}(u)$\index{${\rm
si}(u)$} instead of ${\rm si}(p,u,q)$ and this value is said to be
the {\em rank of semi-isolation}\index{Rank!of semi-isolation!of
label} or the {\em ${\rm si}$-rank}\index{${\rm si}$-rank!of
label} of the label $u$ or of the element $u=\varnothing$ (with
respect to the pair $(p,q)$). For a formula $\theta_{p,u,q}(x,y)$
we set ${\rm si}(\theta_{p,u,q}(x,y))\rightleftharpoons{\rm
si}(u)$.\index{${\rm si}(\theta_{p,u,q}(x,y))$}

\medskip
Clearly, if the theory is small then  the ${\rm si}$-rank of any
label is an ordinal (having a label $u$ with ${\rm
si}(p,u,q)=\infty$, we get continuum many complete types
$r(x,y)\supset p(x)\cup q(y)$).

By the definition we have the following inequality for any formula
$\theta_{p,u,q}(x,y)$ and any realization $a$ of $p$ giving a low
bound for Morley rank of the formula $\theta_{p,u,q}(a,y)$ by the
${\rm si}$-rank:
\begin{equation}\label{si1}
{\rm si}(\theta_{p,u,q}(x,y))\leq{\rm MR}(\theta_{p,u,q}(a,y))+1.
\end{equation}
The inequality (\ref{si1}) implies

\medskip
{\bf Remark 3.1.} If a theory $T$ has a finite Morley rank then
${\rm si}$-ranks of labels $\bigcup\limits_{p,q\in
S^1(T)}\rho_{\nu(p,q)}$ are bounded by the value ${\rm
MR}(x\approx x)+1$.

\medskip
We set ${\rm si}(p,q)\rightleftharpoons{\rm sup}\{{\rm
si}(p,u,q)\mid u\in U\cup\{\varnothing\}\},$\index{${\rm
si}(p,q)$} ${\rm si}(p)\rightleftharpoons{\rm
si}(p,p).$\index{${\rm si}(p)$} For a nonempty family $R$ of
$1$-types, we put ${\rm si}(R)\rightleftharpoons{\rm sup}\{{\rm
si}(p,q)\mid p,q\in R\}$.\index{${\rm si}(R)$} A family $R$ is
called {\em ${\rm si}$-minimal}\index{Family!${\rm si}$-minimal}
if ${\rm si}(R)=1$. The value ${\rm si}(p,q)$ is said to be the
{\em rank of semi-isolation}\index{Rank!of semi-isolation!of type
pair} or the {\em ${\rm si}$-rank}\index{${\rm si}$-rank!of type
pair} of pair $(p,q)$, and ${\rm si}(R)$ is the {\em rank of
semi-isolation}\index{Rank!of semi-isolation!of family of types}
or the {\em ${\rm si}$-rank}\index{${\rm si}$-rank!of family of
types} of the family $R$.

Since there are $|T|$ formulas of a theory $T$ and the inequality
(\ref{si1}) holds we obtain

\medskip
{\bf Proposition 3.2.} {\em Each ${\rm si}$-rank in a theory $T$
is either equal to $\infty$ or less than ${\rm min}\{|T|^+,({\rm
MR}(x\approx x)+1)^+\}$. If Morley rank ${\rm MR}(x\approx x)$ is
equal to an ordinal $\alpha$ then any ${\rm si}$-rank in $T$ is
not more than $\alpha+1$.}

\medskip
The estimation for ${\rm si}$-ranks in Proposition 3.2 can be far
from exact. For instance, ${\rm si}$-ranks in $\omega$-categorical
theories are finite while there are non-$\omega$-stable
$\omega$-categorical theories.

\medskip
{\bf Proposition 3.3.} {\em For any types $p,q\in
S^1(\varnothing)$ the following assertions are satisfied.

$(1)$ If $u,v\in\rho_{\nu(p,q)}\cup\{\varnothing\}$ and $u\unlhd
v$ then ${\rm si}(u)\leq{\rm si}(v)$.

$(2)$ If $u,v\in\rho_{\nu(p,q)}\cup\{\varnothing\}$ then ${\rm
si}(u\vee v)={\rm max}\{{\rm si}(u),{\rm si}(v)\}$ and ${\rm
si}(u\wedge v)\leq{\rm min}\{{\rm si}(u),{\rm si}(v)\}$. The last
inequality is transformed to the equality if and only if there is
a label $v'$ such that $v'\unlhd u$, $v'\unlhd v$, and ${\rm
si}(v')={\rm si}(u)$ or ${\rm si}(v')={\rm si}(v)$.

$(3)$ The equality ${\rm si}(p,q)=0$ holds if and only if there
are no realizations of $p$ semi-isolating realizations of $q$.

$(4)$ The equality ${\rm si}(p,q)=1$ holds if and only if there is
a $(p\rightarrow q)$-formula and each $(p\rightarrow q)$-formula
$\varphi(x,y)$ is equivalent to a disjunction of formulas
$\varphi_i(x,y)$ such that each formula $\varphi_i(a,y)$ is
isolating, where $\models p(a)$.}

\medskip
{\em Proof} is obvious.~$\Box$

\medskip
{\bf Proposition 3.4.} {\em For any nonempty family $R\subseteq
S^1(\varnothing)$ the following assertions are satisfied.

$(1)$ ${\rm si}(R)\geq 1$.

$(2)$ The family $R$ is ${\rm si}$-minimal if and only if for any
types $p,q\in R$ each $(p\rightarrow q)$-formula $\varphi(x,y)$ is
equivalent to a disjunction of formulas $\varphi_i(x,y)$ such that
each formula $\varphi_i(a,y)$ is isolating, where $\models p(a)$.}

\medskip
{\em Proof.} (1) is implied by the inequality ${\rm si}(p)\geq 1$
for any type $p\in S^1(\varnothing)$ since the formula $(a\approx
y)$ witnesses that $a$ semi-isolates itself, where $\models p(a)$.
(2) is an obvious corollary of (1) and Proposition 3.3,
(4).~$\Box$

\medskip
{\bf Remark 3.5.} Since for a strongly minimal theory $T$ the set
of solutions for any formula $\varphi(a,y)$ is finite or cofinite,
any semi-isolating formula $\psi(a,y)$ is represented as a finite
disjunction of some isolating formulas $\psi_i(a,y)$ or as a
negation of a finite disjunction of isolating formulas
$\psi_i(a,y)$. If $\psi(a,y)\vdash p(y)$ and $p(y)$ is a
non-principal type then the representation of $\psi(a,y)$ is
possible only as a finite disjunction of isolating formulas. It
means that ${\rm si}(p)=1$. If $p(y)$ is a principal type and
there are finitely many pairwise non-equivalent isolating formulas
$\psi(a,y)$ with $\models p(a)$ and $\psi(a,y)\vdash p(y)$ then
${\rm si}(p)=1$ too. If there are infinitely many these pairwise
non-equivalent isolating formulas $\psi(a,y)$ then ${\rm
si}(p)=2$.~$\Box$

\medskip
{\bf Definition.} Let $\alpha$ be a positive ordinal, $u_1$ and
$u_2$ be labels in $\rho_{\nu(p,q)}$ such that ${\rm si}(u_1)={\rm
si}(u_2)=\alpha$. The labels $u_1$ and $u_2$ are {\em
$\alpha$-almost identic}\index{Labels!$\alpha$-almost identic} or
{\em
$\sim_\alpha$-equivalent}\index{Labels!$\sim_\alpha$-equivalent}
(denoted by $u_1\sim_\alpha u_2$)\index{$u_1\sim_\alpha u_2$} if
${\rm si}(u_1\div u_2)<\alpha$, where $u_1\div
u_2\rightleftharpoons(u_1\wedge\neg u_2)\vee(\neg u_1\vee
u_2)$.\index{$u_1\div u_2$}

\medskip
{\bf Proposition 3.6.} {\em The relation $\sim_\alpha$ is an
equivalence relation for any set of labels in $\rho_{\nu(p,q)}$
having the ${\rm si}$-rank $\alpha$.}

\medskip
{\em Proof.} Clearly the relation $\sim_\alpha$ is reflexive and
symmetric. For the checking of transitivity we assume that
$u_1\sim_\alpha u_2$ and $u_2\sim_\alpha u_3$. Since
$(u_1\wedge\neg u_2\wedge u_3)\unlhd(u_1\wedge\neg
u_2)\unlhd(u_1\div u_2)$ we have
$${\rm si}(u_1\wedge\neg u_2\wedge u_3)\leq{\rm si}(u_1\div u_2)<\alpha.$$ As $u_1\wedge
u_3=(u_1\wedge u_2\wedge u_3)\vee(u_1\wedge\neg u_2\wedge u_3)$
and ${\rm si}(u_1\wedge\neg u_2\wedge u_3)<\alpha$, for
$u_1\sim_\alpha u_3$, it is enough to prove that ${\rm
si}(u_1\wedge u_2\wedge u_3)=\alpha$. Suppose on contrary that
${\rm si}(u_1\wedge u_2\wedge u_3)<\alpha$. Then ${\rm
si}(u_1\wedge u_2)=\alpha$ and
$$u_1\wedge u_2=(u_1\wedge u_2\wedge u_3)\vee(u_1\wedge
u_2\wedge\neg u_3)$$ imply ${\rm si}(u_1\wedge u_2\wedge\neg
u_3)=\alpha$. But $(u_1\wedge u_2\wedge\neg
u_3)\unlhd(u_2\wedge\neg u_3)\unlhd(u_2\div u_3)$, and ${\rm
si}(u_2\div u_3)<\alpha$ gives ${\rm si}(u_1\wedge u_2\wedge\neg
u_3)<\alpha$. The obtained contradiction means that
$u_1\sim_\alpha u_3$.~$\Box$

\medskip
By the definition, for any label $u\in\rho_{\nu(p,q)}$ having the
${\rm si}$-rank $\alpha$, there is a greatest number
$n\in\omega\setminus\{0\}$ of pairwise inconsistent (or, that
equivalent, of pairwise non-$\sim_\alpha$-equivalent) labels
$u_1,\ldots,u_n$ such that $u_i\unlhd u$ and ${\rm
si}(u_i)=\alpha$, $i=1,\ldots,n$. This number $n$ is called the
{\em degree of semi-isolation}\index{Degree!of semi-isolation!of
label} or the {\em ${\rm si}$-degree}\index{${\rm si}$-degree!of
label} of label $u$ and it is denoted by ${\rm
deg}(p,u,q)$\index{${\rm deg}(p,u,q)$} or by ${\rm
deg}(u)$\index{${\rm deg}(u)$}. We have ${\rm si}(\varnothing)=0$
and put ${\rm deg}(\varnothing)\rightleftharpoons 1$.

\medskip
{\bf Proposition 3.7.} {\em $(1)$ If $u\in\rho_{\nu(p,q)}$ and
${\rm si}(u)=\alpha$ then ${\rm deg}(u)$ is equal to the number of
pairwise inconsistent labels $u_1,\ldots,u_n\in\rho_{\nu(p,q)}$
having the ${\rm si}$-rank $\alpha$, the ${\rm si}$-degree $1$,
and such that $u=u_1\vee\ldots\vee u_n$.

$(2)$ If $u,v\in\rho_{\nu(p,q)}$, ${\rm si}(u)={\rm si}(v)$, and
$u\unlhd v$ then ${\rm deg}(u)\leq{\rm deg}(v)$.

$(3)$ If $u,v\in\rho_{\nu(p,q)}$ and ${\rm si}(u)={\rm si}(v)$
then
$${\rm deg}(u\vee v)\leq{\rm deg}(u)+{\rm deg}(v).$$
The equality in this inequality holds if and only if ${\rm
si}(u\wedge v)<{\rm si}(u)$. If ${\rm si}(u\wedge v)={\rm si}(u)$
then
$${\rm deg}(u\vee v)={\rm deg}(u)+{\rm deg}(v)-{\rm deg}(u\wedge
v).$$

$(4)$ If $u\in\rho_{\nu(p,q)}$ is a label for an isolating
formula, i.~e., $u$ is an atom, then ${\rm si}(u)=1$ and ${\rm
deg}(u)=1$.

$(5)$ If for a label $u\in\rho_{\nu(p,q)}$, ${\rm si}(u)=1$ and
${\rm deg}(u)=1$, then $u$ is not neutral.

$(6).$ If $u\in\rho_{\nu(p,q)}$ and ${\rm si}(u)=1$ then ${\rm
deg}(u)$ is equal to the number of pairwise inconsistent labels
$u_1,\ldots,u_n\in\rho_{\nu(p,q)}$ for isolating formulas such
that $u=u_1\vee\ldots\vee u_n$.}

\medskip
{\em Proof} is obvious.~$\Box$

\medskip
If there is a label $u\in\rho_{\nu(p,q)}$ with ${\rm si}(p,q)={\rm
si}(u)$ then the {\em degree of semi-isolation}\index{Degree!of
semi-isolation!of type pair} or the {\em ${\rm
si}$-degree}\index{${\rm si}$-degree!of type pair} ${\rm
deg}(p,q)$\index{${\rm deg}(p,q)$} of pair $(p,q)$ is
$${\rm sup}\{{\rm deg}(u)\mid u\in\rho_{\nu(p,q)},{\rm
si}(p,q)={\rm si}(u)\},$$ ${\rm deg}(p)\rightleftharpoons{\rm
deg}(p,p)$.

If for a nonempty family $R$ of $1$-types there is a label
$u\in\rho_{\nu(p,q)}$, $p,q\in R$, with ${\rm si}(R)={\rm si}(u)$
then the {\em degree of semi-isolation}\index{Degree!of
semi-isolation!of family of types} or the {\em ${\rm
si}$-degree}\index{${\rm si}$-degree!of family of types} ${\rm
deg}(R)$\index{${\rm deg}(R)$} of $R$ is $${\rm sup}\{{\rm
deg}(u)\mid u\in\rho_{\nu(R)},{\rm si}(R)={\rm si}(u)\}.$$

Clearly, if ${\rm deg}(p,q)$ or ${\rm deg}(R)$ exist then these
values are positive natural numbers or equal $\omega$.

\medskip
For an ordinal $\alpha$, a natural number $n\geq 1$, and a set
$X\in\{U,U\cup\{\varnothing\}\}$ we put
\index{$X\upharpoonright(\alpha,n)$}\index{$X\upharpoonright\alpha$}\index{$X\upharpoonright(\alpha,\omega)$}
$$X\upharpoonright(\alpha,n)\rightleftharpoons\{u\in X\mid{\rm
si}(u)\leq\alpha\mbox{ and if }{\rm si}(u)=\alpha\mbox{ then }{\rm
deg}(u)<n\},$$
$$
X\upharpoonright(\alpha,\omega)\rightleftharpoons
X\upharpoonright\alpha\rightleftharpoons \bigcup\limits_{n\geq
1}X\upharpoonright(\alpha,n).
$$

Clearly, if $\alpha=\beta+1$ then
$X\upharpoonright(\alpha,1)=X\upharpoonright\beta$, and if
$\alpha$ is a limit ordinal then
$X\upharpoonright(\alpha,1)=\bigcup\limits_{\beta<\alpha}X\upharpoonright\beta$.

For ordinals $\alpha,\beta$, $\beta\in(\omega+1)\setminus\{0\}$,
and for the algebra $\mathfrak{A}$ of distributions of binary
semi-isolating formulas of theory $T$ as well as for expansions
and restrictions $\mathfrak{A}'$ of $\mathfrak{A}$, defined in
previous sections, we denote by
$\mathfrak{A}\upharpoonright(\alpha,\beta)$\index{$\mathfrak{A}\upharpoonright(\alpha,\beta)$}
and
$\mathfrak{A}'\upharpoonright(\alpha,\beta)$\index{$\mathfrak{A}'\upharpoonright(\alpha,\beta)$}
as well as by
$\mathfrak{A}_{\alpha,\beta}$\index{$\mathfrak{A}_{\alpha,\beta}$}
and
$\mathfrak{A}'_{\alpha,\beta}$\index{$\mathfrak{A}'_{\alpha,\beta}$}
the restrictions of these algebras to the set
$(U\cup\{\varnothing\})\upharpoonright(\alpha,\beta)$. If
$\beta=\omega$, these restrictions are denoted by
$\mathfrak{A}\upharpoonright\alpha$\index{$\mathfrak{A}\upharpoonright\alpha$},
$\mathfrak{A}'\upharpoonright\alpha$\index{$\mathfrak{A}'\upharpoonright\alpha$},
$\mathfrak{A}_\alpha$\index{$\mathfrak{A}_\alpha$}, and
$\mathfrak{A}'_\alpha$\index{$\mathfrak{A}'_\alpha$}. The
restrictions are called the  {\em
$(\alpha,\beta)$-restrictions}\index{$(\alpha,\beta)$-restriction}
and the {\em $\alpha$-restrictions}\index{$\alpha$-restriction}
respectively.

Since the ${\rm si}$-rank of each label is positive, non-trivial
restrictions (i.~e., with nonempty sets of used labels) are only
the restrictions of algebras with $\alpha>0$. If ${\rm
si}(S^1(\varnothing))=\alpha_0$ then, taking into consideration
the inequality $\alpha>0$, all essential (i.~e., reflecting links
of sets of labels of semi-isolating formulas with respect to their
${\rm si}$-ranks) restrictions of these algebras are formed only
for $0<\alpha<\alpha_0$.

In view of Proposition 3.7 we obtain

\medskip
{\bf Proposition 3.8.} {\em The algebra of distributions of binary
{\sl isolating} formulas of theory $T$ coincides with the algebra
$\mathfrak{A}\upharpoonright(1,2)$. The algebra
$\mathfrak{A}\upharpoonright 1$ consists of labels being
disjunctions of labels of isolating formulas.}

\medskip
\centerline{\bf 4. Monoid of distributions} \centerline{\bf of
binary semi-isolating formulas} \centerline{\bf on a set of
realizations of a type}
\medskip

Consider a complete theory $T$, a type $p(x)\in S(T)$, a regular
labelling function $\nu(p)\mbox{\rm : }{\rm PF}(p)/{\rm PE}(p)\to
U$, and a family of sets ${\rm SI}_p(u_1,\ldots,u_k)$ of labels of
binary semi-isolating formulas, $u_1,\ldots,u_k\in\rho_{\nu(p)}$,
$k\in\omega$.

Below we show some basic properties for sets
$$\lceil u_1,\ldots,u_k\rceil \rightleftharpoons {\rm SI}_p(u_1,\ldots,u_k).$$

\medskip
{\bf Proposition 4.1} {\rm (Associativity).\index{Associativity}}
{\em For any $u_1,u_2,u_3\in\rho_{\nu(p)}$, the following
equalities hold:
$$\lceil\lceil u_1,u_2\rceil,u_3\rceil =\lceil u_1,u_2,u_3\rceil =\lceil u_1,\lceil u_2,u_3\rceil\rceil.$$}

{\em Proof} of inclusions $\lceil\lceil
u_1,u_2\rceil,u_3\rceil\subseteq\lceil u_1,u_2,u_3\rceil$ and
$\lceil u_1,\lceil u_2,u_3\rceil\rceil\subseteq \lceil
u_1,u_2,u_3\rceil$ repeats \cite[Proposition 3.1, 4]{ShS}.

The reverse inclusions are satisfied since, taking labels $v_1$
and $v_2$ for the formulas $\theta_{u_1,u_2}(x,y)$ and
$\theta_{u_2,u_3}(x,y)$, we obtain, for $\models p(a)$, that the
formulas $\theta_{v_1,u_3}(a,y)$, $\theta_{u_1,u_2,u_3}(a,y)$, and
$\theta_{u_1,v_2}(a,y)$ are pairwise equivalent, i.~e.,
$$\lceil v_1,u_3\rceil =\lceil u_1,u_2,u_3\rceil =\lceil u_1,v_2\rceil.\,\,\Box$$

\medskip
In view of associativity, using the induction on number of
parenthesis, we prove that all operations $\lceil
\cdot,\cdot,\ldots,\cdot\rceil$ acting on sets in
$\mathcal{P}(\rho_{\nu(p)})\setminus\{\varnothing\}$ are generated
by the binary operation $\lceil \cdot,\cdot\rceil $ on the set
$\mathcal{P}(\rho_{\nu(p)})\setminus\{\varnothing\}$ and the
values $\lceil X_1,X_2,\ldots,X_k\rceil$,
$X_1,X_2,\ldots,X_k\subseteq\rho_{\nu(p)}$, do not depend on the
sequence of adding of brackets for
$$X_{i,i+1,\ldots,i+m+n}\rightleftharpoons\lceil
X_{i,i+1,\ldots,i+m},X_{i+m+1,i+m+2,\ldots,i+m+n}\rceil,$$ where
$X_{1,2,\ldots,k}=\lceil X_1,X_2,\ldots,X_k\rceil $.

Thus the structure
$\mathfrak{SI}_{\nu(p)}\rightleftharpoons\langle\mathcal{P}(\rho_{\nu(p)})\setminus\{\varnothing\};\lceil
\cdot,\cdot\rceil \rangle$\index{$\mathfrak{SI}_{\nu(p)}$} is a
semigroup admitting the representation of all operations $\lceil
\cdot,\cdot,\ldots,\cdot\rceil$ by terms of the language $\lceil
\cdot,\cdot\rceil$. Further the operation $\lceil
\cdot,\cdot\rceil$ will be denoted also by $\cdot$ and we shall
use the record $uv$ instead of $u\cdot v$.

Since by the choice of the label $0$ for the formula $(x\approx
y)$ the equalities $X\cdot\{0\}=X$ and $\{0\}\cdot X=X$ are true
for any $X\subseteq\rho_{\nu(p)}$, the semigroup
$\mathfrak{SI}_{\nu(p)}$ has the unit $\{0\}$, and it is a monoid.
We have
$$
Y\cdot Z=\bigcup\{yz\mid y\in Y,z\in Z\}
$$
for any sets
$Y,Z\in\mathcal{P}(\rho_{\nu(p)})\setminus\{\varnothing\}$ in this
structure.

Thus the following proposition holds.

\medskip
{\bf Proposition 4.2.} {\em For any complete theory $T$, any type
$p\in S(T)$, and the regular labelling function $\nu(p)$, any
operation ${\rm SI}_p(\cdot,\cdot,\ldots,\cdot)$ on the set
$\mathcal{P}(\rho_{\nu(p)})\setminus\{\varnothing\}$ interpretable
by a term of the monoid $\mathfrak{SI}_{\nu(p)}$.}

\medskip
The monoid $\mathfrak{SI}_{\nu(p)}$ is called the {\em monoid of
binary semi-isolating formulas over the labelling function
$\nu(p)$}\index{Monoid!of binary semi-isolating formulas} or the
{\em ${\rm SI}_{\nu(p)}$-monoid}\index{${\rm
SI}_{\nu(p)}$-monoid}.

\medskip
In view of Propositions 1.3 and 4.1 we obtain

\medskip
{\bf Proposition 4.3.} {\em For any complete theory $T$, any type
$p\in S(T)$, and the regular labelling function $\nu(p)$, the
restriction $\mathfrak{SI}^{\leq
0}_{\nu(p)}$\index{$\mathfrak{SI}^{\leq 0}_{\nu(p)}$} {\rm
(}respectively $\mathfrak{SI}^{\geq 0}_{\nu(p)}$,
$\mathfrak{SI}^{\geq 0,{\rm neu}}_{\nu(p)}${\rm
)}\index{$\mathfrak{SI}^{\geq 0,{\rm neu}}_{\nu(p)}$} of the
monoid $\mathfrak{SI}_{\nu(p)}$ to the set $U^{\leq 0}$ {\rm (}
$U^{\geq 0}$, $U^{\geq 0}\cup U'${\rm )} is a submonoid of
$\mathfrak{SI}_{\nu(p)}$.}

\medskip
By Proposition 3.8, the $(1,2)$-restriction of the monoid
$\mathfrak{SI}_{\nu(p)}$ coincides with the $I_{\nu(p)}$-groupoid
$\mathfrak{P}_{\nu(p)}$. Besides, the $(1,2)$-restrictions of
monoids $\mathfrak{SI}^{\leq 0}_{\nu(p)}$ and $\mathfrak{SI}^{\geq
0}_{\nu(p)}$ equal respectively to the groupoid
$\mathfrak{P}^{\leq 0}_{\nu(p)}$ and the monoid
$\mathfrak{P}^{\geq 0}_{\nu(p)}$.

\medskip
\centerline{\bf 5. $\alpha$-deterministic and almost
$\alpha$-deterministic} \centerline{\bf ${\rm
SI}_{\nu(p)}$-monoids}
\medskip

In the following definition, we generalize the notions of
deterministic and almost deterministic structure
$\mathfrak{P}_{\nu(p)}$ proposed in \cite{ShS}.

\medskip
{\bf Definition.} Let $U_0$ be a subalphabet of the alphabet $U$,
$\alpha$ be a positive ordinal, and $n\geq 1$ be a natural number.
We put
\index{$\rho_{\nu(p),\alpha,n}$}\index{$\rho_{\nu(p),\alpha}$}
$$\rho_{\nu(p),\alpha,n}\rightleftharpoons\{u\in\rho_{\nu(p)}\mid{\rm
si}(u)\leq\alpha,{\rm deg}(u)<n\mbox{ for }{\rm si}(u)=\alpha\},$$
$$\rho_{\nu(p),\alpha}\rightleftharpoons\bigcup\limits_{n\in\omega}\rho_{\nu(p),\alpha,n}.$$

The partial subalgebra $\mathfrak{SI}_{\nu(p)}\upharpoonright U_0$
of the monoid $\mathfrak{SI}_{\nu(p)}$ is called {\em
$(\alpha,n)$-deterministic}\index{Monoid!$(\alpha,n)$-deterministic}\index{Subalgebra!partial!$(\alpha,n)$-deterministic}
if for any labels $u_1,u_2\in\rho_{\nu(p),\alpha,n}\cap U_0$, the
set $\lceil u_1,u_2\rceil\cap U_0$ consists of labels having the
${\rm si}$-ranks $\leq\alpha$ and contains less than $n$ pairwise
non-$\sim_\alpha$-equivalent labels of ${\rm si}$-rank $\alpha$.

The partial subalgebra $\mathfrak{SI}_{\nu(p)}\upharpoonright U_0$
of the monoid $\mathfrak{SI}_{\nu(p)}$ is called {\em
$\alpha$-deterministic}\index{Monoid!$\alpha$-deterministic}\index{Subalgebra!partial!$\alpha$-deterministic}
if $\mathfrak{SI}_{\nu(p)}\upharpoonright U_0$ is
$(\alpha,2)$-deterministic.

The partial subalgebra $\mathfrak{SI}_{\nu(p)}\upharpoonright U_0$
of the monoid $\mathfrak{SI}_{\nu(p)}$ is called {\em almost
$\alpha$-deterministic}\index{Monoid!almost
$\alpha$-deterministic}\index{Subalgebra!partial!almost
$\alpha$-deterministic} or {\em
$(\alpha,\omega)$-deterministic}\index{Monoid!$(\alpha,\omega)$-deterministic}\index{Subalgebra!partial!$(\alpha,\omega)$-deterministic}
if for any labels $u_1,u_2\in\rho_{\nu(p),\alpha}\cap U_0$, the
set $\lceil u_1,u_2\rceil\cap U_0$ consists of labels having the
${\rm si}$-ranks $\leq\alpha$ and contains finitely many pairwise
non-$\sim_\alpha$-equivalent labels of ${\rm si}$-rank $\alpha$.

\medskip
By the definition, each $(\alpha,\omega)$-deterministic structure
$\mathfrak{SI}_{\nu(p)}\upharpoonright U_0$ is a union of its
$(\alpha,n)$-deterministic substructures, $n\geq 1$. So each
$\alpha$-deterministic structure
$\mathfrak{SI}_{\nu(p)}\upharpoonright U_0$ is almost
$\alpha$-deterministic.

If $U_0=U$ we shall not point out restrictions to the set $U_0$
for considered structures.

Below we show some basic properties of (almost)
$\alpha$-deterministic partial algebras
$\mathfrak{SI}_{\nu(p)}\upharpoonright U_0$.

\medskip
{\bf Proposition 5.1.} {\em {\rm (Monotony)\index{Monotony}} If a
structure $\mathfrak{SI}_{\nu(p)}\upharpoonright U_0$ is {\rm
(}almost{\rm )} $\alpha$-deterministic and $\beta$ is a positive
ordinal then the structure $(\mathfrak{SI}_{\nu(p)}\upharpoonright
U_0)\upharpoonright\beta$ is also {\rm (}almost{\rm )}
$\alpha$-deterministic.}

\medskip
{\em Proof} is obvious.~$\Box$

\medskip
{\bf Proposition 5.2.} {\em For any monoid
$\mathfrak{SI}_{\nu(p)}$ and ordinals $\alpha,\beta$, where
$\alpha,\beta>0$, $\beta\in\omega+1$, the following conditions are
equivalent:

{\rm (1)} the monoid $\mathfrak{SI}_{\nu(p)}$ is
$(\alpha,\beta)$-deterministic;

{\rm (2)} ${\rm si}(u_1\circ u_2)\leq\alpha$ for any labels
$u_1,u_2\in\rho_{\nu(p),\alpha,\beta}$ and if ${\rm si}(u_1\circ
u_2)=\alpha$ then ${\rm deg}(u_1\circ u_2)<\beta$.}

\medskip
{\em Proof.} The implication $(1)\Rightarrow(2)$ is obvious.

$(2)\Rightarrow(1)$. Consider arbitrary labels
$u_1,u_2\in\rho_{\nu(p),\alpha,\beta}$. Since, by hypothesis,
${\rm si}(u_1\circ u_2)\leq\alpha$ and $v\unlhd(u_1\circ u_2)$ for
any label $v\in\lceil u_1,u_2\rceil$, $\lceil u_1,u_2\rceil$
consists of labels of ${\rm si}$-ranks $\leq\alpha$, and if ${\rm
si}(v)={\rm si}(u_1\circ u_2)=\alpha$ then ${\rm deg}(v)\leq{\rm
deg}(u_1\circ u_2)<\beta$. Thus, the monoid
$\mathfrak{SI}_{\nu(p)}$ is $(\alpha,\beta)$-deterministic.~$\Box$

\medskip
Proposition 5.2 immediately implies

\medskip
{\bf Corollary 5.3.} {\em For any monoid $\mathfrak{SI}_{\nu(p)}$
and a positive ordinal $\alpha$ the following conditions are
equivalent:

{\rm (1)} the monoid $\mathfrak{SI}_{\nu(p)}$ is almost
$\alpha$-deterministic;

{\rm (2)} ${\rm si}(u_1\circ u_2)\leq\alpha$ for any labels
$u_1,u_2\in\rho_{\nu(p),\alpha}$.}

\medskip
{\bf Corollary 5.4.} {\em If ${\rm si}(p)$ is an ordinal then the
monoid $\mathfrak{SI}_{\nu(p)}$ is almost ${\rm
si}(p)$-deterministic.}

\medskip
{\bf Proposition 5.5.} {\em If a monoid $\mathfrak{SI}_{\nu(p)}$
is $(\alpha,\beta)$-deterministic then the structure
$\mathfrak{SI}_{\nu(p),\alpha}\rightleftharpoons
\mathfrak{SI}_{\nu(p)}\upharpoonright\alpha$\index{$\mathfrak{SI}_{\nu(p),\alpha}$}
is also an $(\alpha,\beta)$-deterministic monoid.}

\medskip
{\em Proof.} Since for any $\alpha$-restriction the associativity,
the presence of unit $\{0\}$, and the $(\alpha,\beta)$-determinacy
is preserved, it is enough to note that for any labels $u_1$ and
$u_2$ in $\mathfrak{SI}_{\nu(p),\alpha,\beta}$ there is a label
$v$ in $\mathfrak{SI}_{\nu(p),\alpha,\beta}$ belonging $\lceil
u_1,u_2\rceil$. We can take $u_1\circ u_2$ for $v$ since, by
hypothesis, ${\rm si}(u_1\circ u_2)\leq\alpha$ and if ${\rm
si}(u_1\circ u_2)=\alpha$ then ${\rm deg}(u_1\circ
u_2)<\beta$.~$\Box$

\medskip
{\bf Proposition 5.6.} {\em If ${\rm si}(p)$ is an ordinal then
the monoid $\mathfrak{SI}_{\nu(p)}$ is ${\rm si}(p)$-deterministic
if and only if the value ${\rm deg}(p)$ is not defined or equals
$1$.}

\medskip
{\em Proof.} If ${\rm deg}(p)$ is not defined the ordinal
$\alpha={\rm si}(p)$ is limit and can not be achieved by labels in
$\rho_{\nu(p)}$. In particular, for any $u_1,u_2\in\rho_{\nu(p)}$
the set $\lceil u_1,u_2\rceil$ does not contain labels of ${\rm
si}$-rank $\alpha$. If ${\rm deg}(p)\geq 2$ then there are
non-$\sim_\alpha$-equivalent labels $u_1,u_2\in\rho_{\nu(p)}$ of
${\rm si}$-rank $\alpha$. Then $\lceil u_1\vee u_2,0\rceil$
contains the labels $u_1$ and $u_2$, whence the monoid
$\mathfrak{SI}_{\nu(p)}$ is not $\alpha$-deterministic. If ${\rm
deg}(p)=1$ then there is unique, up to $\sim_\alpha$-equivalence,
label in $\rho_{\nu(p)}$ having the ${\rm si}$-rank $\alpha$.
Since such a label is unique, the monoid $\mathfrak{SI}_{\nu(p)}$
is $\alpha$-deterministic.~$\Box$

\medskip
{\bf Proposition 5.7.} {\em The structure $\mathfrak{P}_{\nu(p)}$
is {\rm (}almost{\rm )} deterministic if and only if the structure
$\mathfrak{SI}_{\nu(p),1,2}$ is {\rm (}almost{\rm )}
$1$-deterministic.}

\medskip
{\em Proof} follows by the equality
$\mathfrak{SI}_{\nu(p),1,2}=\mathfrak{P}_{\nu(p)}$.~$\Box$

\medskip
{\bf Proposition 5.8.} {\em Let $p(x)$ be a complete type of a
theory $T$, $\nu(p)$ be a regular labelling function, and ${\rm
si}(p)<\omega$. The following conditions are equivalent:

$(1)$ the monoid $\mathfrak{SI}_{\nu(p)}$ is $(1,n)$-deterministic
for some $n\in\omega$;

$(2)$ the set $\rho_{\nu(p)}$ is finite;

$(3)$ the set $\rho_{\nu(p),1}$ finite;

$(4)$ the set $\rho_{\nu(p),1,2}$ {\rm (}consisting of all atoms
$u\in\rho_{\nu(p)}${\rm )} is finite.}

\medskip
{\em Proof.} If ${\rm si}(p)>1$ then, by ${\rm si}(p)<\omega$, the
set $\rho_{\nu(p),1}$ is infinite and so the set $\rho_{\nu(p)}$
is also infinite. Since each label in $\rho_{\nu(p),1}$ is a
disjunction of labels in $\rho_{\nu(p),1,2}$ and for any labels
$u_1,\ldots,u_n\in\rho_{\nu(p),1}$ the label $u_1\vee\ldots\vee
u_n$ belongs also to $\rho_{\nu(p),1}$, the set
$\rho_{\nu(p),1,2}$ is infinite and the monoid
$\mathfrak{SI}_{\nu(p)}$ is not $(1,n)$-deterministic for
$n\in\omega$. Thus, none of the conditions (1)--(4) is not
satisfied.

If ${\rm si}(p)=1$ then each label in $\rho_{\nu(p)}$ has the
${\rm si}$-rank $1$ and is represented as a disjunction of labels
in $\rho_{\nu(p),1,2}$. Thus, the conditions (2)--(4) are
equivalent. If the set $\rho_{\nu(p),1,2}$ contains $m\in\omega$
labels then there are $2^m-1$ labels forming the set
$\rho_{\nu(p)}$. Hence, the monoid $\mathfrak{SI}_{\nu(p)}$ is
$(1,2^m-1)$-deterministic. If the set $\rho_{\nu(p),1,2}$ is
infinite then, for pairwise distinct labels
$u_1,\ldots,u_m\in\rho_{\nu(p),1,2}$, the set $\lceil
u_1\vee\ldots u_m,0\rceil$ contains $2^m-1$ labels and, since $m$
is not bounded, the monoid $\mathfrak{SI}_{\nu(p)}$ is not
$(1,n)$-deterministic for any $n$. Thus, the condition (1) is
equivalent to each of the conditions $(2)$--$(4)$.~$\Box$

\medskip
Proposition 5.8 and \cite[Corollary 7.4]{ShS} imply

\medskip
{\bf Corollary 5.9.} {\em Let $p(x)$ be a complete type of a
theory $T$, $\nu(p)$ be a regular labelling function, ${\rm
si}(p)<\omega$, and $\mathfrak{SI}_{\nu(p)}$ is a
$(1,n)$-deterministic monoid, for some $n\in\omega$, having a
negative label. Then the groupoid $\mathfrak{SI}_{\nu(p),1,2}$
generates the strict order property.}

\medskip
{\bf Definition} \cite{SuLP, Su92}. Let $p(x)$ be a type in
$S(T)$. A type $q(x_1,\ldots,x_n)\in S(T)$ is called a
\emph{$(n,p)$-type}\index{$(n,p)$-type} if
$q(x_1,\ldots,x_n)\supseteq\bigcup\limits_{i=1}^n p(x_i)$. The set
of all $(n,p)$-types of $T$ is denoted by
$S_{n,p}(T)$\index{$S_{n,p}(T)$} and elements of the set
$S_p(T)\rightleftharpoons\bigcup\limits_{n\in\omega\setminus\{0\}}S_{n,p}(T)$\index{$S_p(T)$}
are called \emph{$p$-types}\index{$p$-type}.

A type $q(\bar{y})$ in $S_p(T)$ is called
\emph{$p$-principal}\index{Type!$p$-principal} if there is a
formula $\varphi(\bar{y})\in q(\bar{y})$ such that
$\cup\{p(y_i)\mid y_i\in\bar{y}\}\cup\{\varphi(\bar{y})\}\vdash
q(\bar{y})$.

\medskip
{\bf Lemma 5.10} {\rm \cite{SuLP, Su92}.} {\em For any type $p$
and a natural number $n\geq 1$ the following conditions are
equivalent:

{\rm (1)} the set of $(n,p)$-types with a tuple $(x_1,\ldots,x_n)$
of free variables is infinite;

{\rm (2)} there is a non-$p$-principal $(n,p)$-type.}

\medskip
Proposition 5.8 and Lemma 5.10 imply

\medskip
{\bf Corollary 5.11.} {\em If $p(x)$ is a complete type of a
theory $T$, $\nu(p)$ is a regular labelling function, and all
$(2,p)$-types are $p$-principal, then the monoid
$\mathfrak{SI}_{\nu(p)}$ is $(1,n)$-deterministic for some
$n\in\omega$.}

\medskip
By Corollaries 5.9 and 5.11, we obtain

\medskip
{\bf Corollary 5.12.} {\em Let $p(x)$ be a complete type of a
theory $T$, $\nu(p)$ be a regular labelling function,
$\rho_{\nu(p)}\cap U^-\ne\varnothing$, and all $(2,p)$-types be
$p$-principal. Then the groupoid $\mathfrak{SI}_{\nu(p),1,2}$
generates the strict order property.}

\medskip
For a type $p(x)$ and a positive ordinal $\alpha$, we denote by
${\rm SI}_{p,\alpha}$\index{${\rm SI}_{p,\alpha}$} (in a model
$\mathcal{M}$ of $T$) the relation of semi-isolation (over
$\varnothing$) on a set of realizations of $p$ restricted to the
set of formulas of ${\rm si}$-rank $\leq\alpha$:
$$
{\rm
SI}_{p,\alpha}\rightleftharpoons\{(a,b)\mid\:\mathcal{M}\models
p(a)\wedge p(b)\mbox{ and }a\mbox{ semi-isolates }b$$
$$\mbox{ by a formula }\theta_u(x,y)\mbox{ with a}{\rm
si}\mbox{-rank }\leq\alpha\}.
$$

Clearly, $I_p={\rm SI}_{p,1}$ for any type $p\in
S^1(\varnothing)$. Seeing this equality and
$\mathfrak{SI}_{\nu(p),1,2}=\mathfrak{P}_{\nu(p)}$ the following
proposition generalizes Proposition 4.3 in \cite{ShS}.

\medskip
{\bf Proposition 5.13.} {\em Let $p(x)$ be a complete type of a
theory $T$, $\nu(p)$ be a regular labelling function, and $\alpha$
be a positive ordinal. The following conditions are equivalent:

$(1)$ the relation ${\rm SI}_{p,\alpha}$ {\rm (}on a set of
realizations of $p$ in any model $\mathcal{M}\models T${\rm )} is
transitive;

$(2)$ the structure $\mathfrak{SI}_{\nu(p),\alpha}$ is an almost
$\alpha$-deterministic monoid.}

\medskip
{\em Proof.} Let $a$, $b$, and $c$ be realizations of $p$ such
that $(a,b)\in {\rm SI}_{p,\alpha}$ and $(b,c)\in {\rm
SI}_{p,\alpha}$ by semi-isolating formulas $\theta_{u_1}(a,y)$ and
$\theta_{u_2}(b,y)$. If the structure
$\mathfrak{SI}_{\nu(p),\alpha}$ is an almost
$\alpha$-deterministic monoid then ${\rm si}(u_1\circ
u_2)\leq\alpha$ and the pair $(a,c)$ belongs to ${\rm
SI}_{p,\alpha}$ by the semi-isolating formula
$\theta_{u_1,u_2}(x,y)$. Since elements $a$ $,b$, and $c$ are
arbitrary we have $(2)\Rightarrow(1)$.

Assume now that for some $u_1,u_2\in\rho_{\nu(p),\alpha}$ the set
${\rm SI}_p(u_1,u_2)$ contains a label $u$ such that ${\rm
si}(u)>\alpha$. Then by compactness the set
$$q(a,y)\rightleftharpoons\{\theta_{u_1,u_2}(a,y)\}\cup\{\neg\theta_v(a,y)\mid
v\in{\rm SI}_p(u_1,u_2),{\rm si}(v)\leq\alpha\}$$ is consistent,
where $\models p(a)$. Consider realizations $b$ and $c$ of $p$
such that $\models\theta_{u_1}(a,b)\wedge\theta_{u_2}(b,c)$ and
$\models q(a,c)$. We have $(a,b)\in {\rm SI}_{p,\alpha}$ and
$(b,c)\in {\rm SI}_{p,\alpha}$ but $(a,c)\notin {\rm
SI}_{p,\alpha}$ by the construction of $q$. Thus, the relation
${\rm SI}_{p,\alpha}$ is not transitive and the implication
$(1)\Rightarrow(2)$ holds.~$\Box$

\medskip
Note that for any ordinal $\alpha>0$ there are no
$(p,\theta_u,p)$-edges, linking distinct realizations of $p$ and
satisfying the conditions $u>0$, ${\rm si}(u)\leq\alpha$, and
${\rm si}(u^{-1})\leq\alpha$, if and only if the relation ${\rm
SI}_{p,\alpha}$ is antisymmetric. Since ${\rm SI}_{p,\alpha}$ is
reflexive, the definition of $\nu(p)$ and Propositions 1.3, 5.13
imply

\medskip
{\bf Corollary 5.14.} {\em Let $p(x)$ be a complete type of a
theory $T$, $\nu(p)$ be a regular labelling function, and $\alpha$
be a positive ordinal. The following conditions are equivalent:

$(1)$ the relation ${\rm SI}_{p,\alpha}$ is a partial order on a
set of realizations of $p$ in any model $\mathcal{M}\models T$;

$(2)$ the structure $\mathfrak{SI}_{\nu(p),\alpha}$ is an almost
$\alpha$-deterministic monoid and $\rho_{\nu(p),\alpha}\subseteq
U^{\leq 0}$.

This partial order $\mathfrak{SI}_{\nu(p),\alpha}$ is identical if
and only if $\rho_{\nu(p),\alpha}=\{0\}$. If ${\rm SI}_{p,\alpha}$
is not identical, it has infinite chains.}

\medskip
Propositions 1.3 and 5.13 imply also

\medskip
{\bf Corollary 5.15.} {\em Let $p(x)$ be a complete type of a
theory $T$, $\nu(p)$ be a regular labelling function, and $\alpha$
be a positive ordinal. The following conditions are equivalent:

$(1)$ the relation ${\rm SI}_{p,\alpha}$ is an equivalence
relation on the set of realizations of $p$ in any model
$\mathcal{M}\models T$;

$(2)$ the structure $\mathfrak{SI}_{\nu(p),\alpha}$ is an almost
$\alpha$-deterministic monoid and consists of labels in $U^{\geq
0}$.}

\medskip
Recall \cite{ShS} that an element $u\in \rho_{\nu(p)}$ is called
{\em {\rm (}almost{\em )}
deterministic}\index{Element!deterministic}\index{Element!almost
deterministic} if for any/some realization $a$ of $p$ the formula
$\theta_u(a,y)$ has unique solution (has finitely many solutions).

\medskip
Since each semi-isolating formula $\theta_u(a,y)$ with finitely
many solutions is equivalent to a disjunction of isolating
formulas $\theta_{u_i}(a,y)$, each almost deterministic element
has the ${\rm si}$-rank $1$ and so belongs to the set of labels in
the structure $\mathfrak{SI}_{\nu(p),1,n+1}$, where $n$ is the
number of solutions for $\theta_u(a,y)$. In particular, each
deterministic element belongs to the set of labels in the
structure $\mathfrak{SI}_{\nu(p),1,2}$.

It is shown in \cite[Proposition 4.7]{ShS} that if elements $u$
and $v$ are {\rm (}almost{\rm )} deterministic then each element
$v'$ in $u\cdot v$ is also {\rm (}almost{\rm )} deterministic.
Hence, the ${\rm si}$-rank $1$ is preserved for compositions
$u\circ v$ of {\rm (}almost{\rm )} deterministic elements $u$ and
$v$. Moreover, the ${\rm si}$-degree $1$ is preserved for
compositions of deterministic elements.

\begin{figure}[t]
\begin{center}
\unitlength 1.7cm
\begin{picture}(8,6)(1.5,-2.0)
{\footnotesize

\put(1.8,-1){\makebox(0,0)[cc]{$\bullet$}}
\put(1.8,0){\makebox(0,0)[cc]{$\bullet$}}
\put(1.8,1){\makebox(0,0)[cc]{$\bullet$}}
\put(1.8,2){\makebox(0,0)[cc]{$\bullet$}}
\put(1.8,2.2){\makebox(0,0)[cc]{$\cdot$}}
\put(1.8,2.4){\makebox(0,0)[cc]{$\cdot$}}
\put(1.8,2.6){\makebox(0,0)[cc]{$\cdot$}}
\put(1.8,2.8){\makebox(0,0)[cc]{$\cdot$}}
\put(1.8,3){\makebox(0,0)[cc]{$\cdot$}}
\put(1.8,3.2){\makebox(0,0)[cc]{$\cdot$}}
\put(1.8,3.4){\makebox(0,0)[cc]{$\cdot$}}
\put(1.8,3.6){\makebox(0,0)[cc]{$\bullet$}}
\put(1.8,-1){\line(0,1){2}} 
\put(1.8,1){\line(0,1){1}}
\put(1.9,-1.1){\makebox(0,0)[cl]{$\{0\}$}}
\put(1.9,0){\makebox(0,0)[cl]{$\mathfrak{SI}_{1}$}}
\put(1.9,1){\makebox(0,0)[cl]{$\mathfrak{SI}_{2}$}}
\put(1.9,2){\makebox(0,0)[cl]{$\mathfrak{SI}_{3}$}}
\put(1.9,3.7){\makebox(0,0)[cl]{$\mathfrak{SI}$}}
\put(1.8,-1.6){\makebox(0,0)[cc]{\small\em a}}

\put(3,0){\makebox(0,0)[cc]{$\bullet$}}
\put(3,1){\makebox(0,0)[cc]{$\bullet$}}
\put(3,1.2){\makebox(0,0)[cc]{$\cdot$}}
\put(3,1.4){\makebox(0,0)[cc]{$\cdot$}}
\put(3,1.6){\makebox(0,0)[cc]{$\cdot$}}
\put(3,1.8){\makebox(0,0)[cc]{$\cdot$}}
\put(3,2){\makebox(0,0)[cc]{$\cdot$}}
\put(3,2.2){\makebox(0,0)[cc]{$\cdot$}}
\put(3,2.4){\makebox(0,0)[cc]{$\cdot$}}
\put(3,2.6){\makebox(0,0)[cc]{$\bullet$}}
\put(3,3.6){\makebox(0,0)[cc]{$\bullet$}}
\put(4,-1){\makebox(0,0)[cc]{$\bullet$}}
\put(4,1){\makebox(0,0)[cc]{$\bullet$}}
\put(4,2){\makebox(0,0)[cc]{$\bullet$}}
\put(4,2.2){\makebox(0,0)[cc]{$\cdot$}}
\put(4,2.4){\makebox(0,0)[cc]{$\cdot$}}
\put(4,2.6){\makebox(0,0)[cc]{$\cdot$}}
\put(4,2.8){\makebox(0,0)[cc]{$\cdot$}}
\put(4,3){\makebox(0,0)[cc]{$\cdot$}}
\put(4,3.2){\makebox(0,0)[cc]{$\cdot$}}
\put(4,3.4){\makebox(0,0)[cc]{$\cdot$}}
\put(4,3.6){\makebox(0,0)[cc]{$\bullet$}}
\put(5,0){\makebox(0,0)[cc]{$\bullet$}}
\put(5,1){\makebox(0,0)[cc]{$\bullet$}}
\put(5,1.2){\makebox(0,0)[cc]{$\cdot$}}
\put(5,1.4){\makebox(0,0)[cc]{$\cdot$}}
\put(5,1.6){\makebox(0,0)[cc]{$\cdot$}}
\put(5,1.8){\makebox(0,0)[cc]{$\cdot$}}
\put(5,2){\makebox(0,0)[cc]{$\cdot$}}
\put(5,2.2){\makebox(0,0)[cc]{$\cdot$}}
\put(5,2.4){\makebox(0,0)[cc]{$\cdot$}}
\put(5,2.6){\makebox(0,0)[cc]{$\bullet$}}
\put(5,3.6){\makebox(0,0)[cc]{$\bullet$}}
\put(3,2.6){\line(0,1){1}} \put(3,1){\line(0,-1){1}}
\put(3,2.6){\line(1,1){1}} \put(5,1){\line(0,-1){0.5}}
\put(5,2.6){\line(0,1){1}}
\put(4,0){\line(0,1){1}} 
\put(5,1){\line(0,-1){1}} \put(4,1){\line(0,-1){1}}
\put(4,1){\line(0,1){1}} \put(4,0){\line(0,-1){1}}
\put(4,1){\line(-1,-1){1}} \put(4,1){\line(1,-1){1}}
\put(4,2){\line(-1,-1){1}} \put(4,3.6){\line(1,-1){1}}
\put(4,2){\line(1,-1){1}} \put(3,0){\line(1,-1){1}}
\put(5,0){\line(-1,-1){1}}
\put(3.13,0){\makebox(0,0)[cl]{$\mathfrak{SI}^{\leq 0}_{1,2}$}}
\put(4.1,-1.1){\makebox(0,0)[cl]{$\{0\}$}}
\put(5.1,0){\makebox(0,0)[cl]{$\mathfrak{SI}^{\geq 0}_{1,2}$}}
\put(3.09,1){\makebox(0,0)[cl]{$\mathfrak{SI}^{\leq 0}_{1,3}$}}
\put(4.1,1){\makebox(0,0)[cl]{$\mathfrak{SI}_{1,2}$}}
\put(5.1,1){\makebox(0,0)[cl]{$\mathfrak{SI}^{\geq 0}_{1,3}$}}
\put(3.13,2.6){\makebox(0,0)[cl]{$\mathfrak{SI}^{\leq 0}_{1}$}}
\put(4.1,2){\makebox(0,0)[cl]{$\mathfrak{SI}_{1,3}$}}
\put(5.1,2.6){\makebox(0,0)[cl]{$\mathfrak{SI}^{\geq 0}_{1}$}}
\put(4.1,3.7){\makebox(0,0)[cl]{$\mathfrak{SI}_{1}$}}
\put(3.1,3.7){\makebox(0,0)[cl]{$\mathfrak{SI}^{\leq 0}_{2,2}$}}
\put(5.1,3.7){\makebox(0,0)[cl]{$\mathfrak{SI}^{\geq 0}_{2,2}$}}

\put(4,-1.6){\makebox(0,0)[cc]{\small\em b}}

\put(6.2,-1){\makebox(0,0)[cc]{$\bullet$}}
\put(6.2,0){\makebox(0,0)[cc]{$\bullet$}}
\put(6.2,1){\makebox(0,0)[cc]{$\bullet$}}
\put(6.2,1.2){\makebox(0,0)[cc]{$\cdot$}}
\put(6.2,1.4){\makebox(0,0)[cc]{$\cdot$}}
\put(6.2,1.6){\makebox(0,0)[cc]{$\cdot$}}
\put(6.2,1.8){\makebox(0,0)[cc]{$\cdot$}}
\put(6.2,2){\makebox(0,0)[cc]{$\cdot$}}
\put(6.2,2.2){\makebox(0,0)[cc]{$\cdot$}}
\put(6.2,2.4){\makebox(0,0)[cc]{$\cdot$}}
\put(6.2,2.6){\makebox(0,0)[cc]{$\bullet$}}
\put(6.2,3.6){\makebox(0,0)[cc]{$\bullet$}}
\put(7.2,-1){\makebox(0,0)[cc]{$\bullet$}}
\put(7.2,0){\makebox(0,0)[cc]{$\bullet$}}
\put(7.2,1){\makebox(0,0)[cc]{$\bullet$}}
\put(7.2,2){\makebox(0,0)[cc]{$\bullet$}}
\put(7.2,2.2){\makebox(0,0)[cc]{$\cdot$}}
\put(7.2,2.4){\makebox(0,0)[cc]{$\cdot$}}
\put(7.2,2.6){\makebox(0,0)[cc]{$\cdot$}}
\put(7.2,2.8){\makebox(0,0)[cc]{$\cdot$}}
\put(7.2,3){\makebox(0,0)[cc]{$\cdot$}}
\put(7.2,3.2){\makebox(0,0)[cc]{$\cdot$}}
\put(7.2,3.4){\makebox(0,0)[cc]{$\cdot$}}
\put(7.2,3.6){\makebox(0,0)[cc]{$\bullet$}}
\put(8.2,-1){\makebox(0,0)[cc]{$\bullet$}}
\put(8.2,0){\makebox(0,0)[cc]{$\bullet$}}
\put(8.2,1){\makebox(0,0)[cc]{$\bullet$}}
\put(8.2,1.2){\makebox(0,0)[cc]{$\cdot$}}
\put(8.2,1.4){\makebox(0,0)[cc]{$\cdot$}}
\put(8.2,1.6){\makebox(0,0)[cc]{$\cdot$}}
\put(8.2,1.8){\makebox(0,0)[cc]{$\cdot$}}
\put(8.2,2){\makebox(0,0)[cc]{$\cdot$}}
\put(8.2,2.2){\makebox(0,0)[cc]{$\cdot$}}
\put(8.2,2.4){\makebox(0,0)[cc]{$\cdot$}}
\put(8.2,2.6){\makebox(0,0)[cc]{$\bullet$}}
\put(8.2,3.6){\makebox(0,0)[cc]{$\bullet$}}
\put(6.2,2.6){\line(0,1){1}} \put(6.2,1){\line(0,-1){1}}
\put(6.2,2.6){\line(1,1){1}} \put(8.2,1){\line(0,-1){0.5}}
\put(8.2,2.6){\line(0,1){1}}
\put(7.2,0){\line(0,1){1}} 
\put(7.2,0){\line(-1,-1){1}} \put(7.2,0){\line(1,-1){1}}
\put(8.2,1){\line(0,-1){1}} \put(7.2,1){\line(0,-1){1}}
\put(7.2,1){\line(0,1){1}} \put(7.2,0){\line(0,-1){1}}
\put(7.2,1){\line(-1,-1){1}} \put(7.2,1){\line(1,-1){1}}
\put(7.2,2){\line(-1,-1){1}} \put(7.2,3.6){\line(1,-1){1}}
\put(7.2,2){\line(1,-1){1}} \put(6.2,0){\line(0,-1){1}}
\put(8.2,0){\line(0,-1){1}}

\put(7.3,-1.1){\makebox(0,0)[cl]{$\mathfrak{SI}_{\alpha}$}}
\put(6.3,-1.1){\makebox(0,0)[cl]{$\mathfrak{SI}^{\leq
0}_{\alpha+1,2}$}}
\put(8.3,-1.1){\makebox(0,0)[cl]{$\mathfrak{SI}^{\geq
0}_{\alpha+1,2}$}}
\put(6.32,0){\makebox(0,0)[cl]{$\mathfrak{SI}^{\leq
0}_{\alpha+_1,3}$}}
\put(7.32,0){\makebox(0,0)[cl]{$\mathfrak{SI}_{\alpha+1,2}$}}
\put(8.3,0){\makebox(0,0)[cl]{$\mathfrak{SI}^{\geq
0}_{\alpha+1,3}$}}
\put(6.29,1){\makebox(0,0)[cl]{$\mathfrak{SI}^{\leq
0}_{\alpha+1,4}$}}
\put(7.3,1){\makebox(0,0)[cl]{$\mathfrak{SI}_{\alpha+1,3}$}}
\put(8.3,1){\makebox(0,0)[cl]{$\mathfrak{SI}^{\geq
0}_{\alpha+1,4}$}}
\put(6.33,2.6){\makebox(0,0)[cl]{$\mathfrak{SI}^{\leq
0}_{\alpha+1}$}}
\put(7.3,2){\makebox(0,0)[cl]{$\mathfrak{SI}_{\alpha+1,4}$}}
\put(8.3,2.6){\makebox(0,0)[cl]{$\mathfrak{SI}^{\geq
0}_{\alpha+1}$}}
\put(7.3,3.7){\makebox(0,0)[cl]{$\mathfrak{SI}_{\alpha+1}$}}
\put(6.3,3.7){\makebox(0,0)[cl]{$\mathfrak{SI}^{\leq
0}_{\alpha+2,2}$}}
\put(8.3,3.7){\makebox(0,0)[cl]{$\mathfrak{SI}^{\geq
0}_{\alpha+2,2}$}}

\put(7.2,-1.6){\makebox(0,0)[cc]{\small\em c}} }
\end{picture}
\vspace{1cm} Fig. 1
\end{center}
\end{figure}

In Figure 1, the fragments of Hasse diagram are presented
illustrating the links of the structure
$\mathfrak{SI}\rightleftharpoons\mathfrak{SI}_{\nu(p)}$ with
structures above, being restrictions of $\mathfrak{SI}$ to
subalphabets of $U$. Here the superscripts $\cdot^{\leq 0}$ and
$\cdot^{\geq 0}$ point out on restrictions of $\mathfrak{SI}$ to
the sets $U^{\leq 0}$ and $U^{\geq 0}$ respectively, and the
subscripts to the upper estimates for ${\rm si}$-ranks and ${\rm
si}$-degrees of labels. In Figure~1,~a, a hierarchy of structures
$\mathfrak{SI}_{\alpha}$, $\alpha\leq{\rm si}(p)$, is depicted
starting with the trivial substructure; in Figure~1,~b, links
between substructures of $\mathfrak{SI}_{\nu(p),1}$ are presented;
in Figure~1,~c, links between substructures of
$\mathfrak{SI}_{\alpha+1}$ for $1\leq\alpha<{\rm si}(p)$ are
shown. For a limit ordinal $\beta\leq{\rm si}(p)$, the Hasse
diagram for substructures of $\mathfrak{SI}_\beta$ is obtained by
union of presented diagrams for $\alpha<\beta$. If an ordinal
$\beta\leq{\rm si}(p)$ is not limit, the Hasse diagram corresponds
to the union of presented diagrams for $\alpha<\beta$ with the
removal of structures $\mathfrak{SI}^{\leq 0}_{\beta+1,2}$ and
$\mathfrak{SI}^{\geq 0}_{\beta+1,2}$.

\medskip
\centerline{\bf 6. ${\rm POSTC}$-monoids}
\medskip

In this Section, we shall consider both the monoids
$\mathfrak{SI}_{\nu(p)}$ and their expansions (with the addition
of empty set to the universe such that
$X\cdot\varnothing=\varnothing\cdot X=\varnothing$ for
$X\in\mathcal{P}(\rho_{\nu(p)})$ \footnote{This extension forms
also a monoid with the unit $\varnothing$ instead of $\{0\}$.}) by
operations and relations of ${\rm POSTC}$-algebras containing
these monoids. These expansions\index{$\mathfrak{M}_{\nu(p)}$}
$$
\mathfrak{M}_{\nu(p)}\rightleftharpoons\langle\mathcal{P}(\rho_{\nu(p)});\,\cdot,\unlhd,\vee,\wedge,(\cdot\,\wedge\,\neg\,
\cdot),\circ\rangle
$$
are called {\em preordered monoids with relative set-theoretic
operations and compositions}\index{Monoid!preordered with relative
set-theoretic operations and a composition} over regular labelling
functions $\nu(p)$, or briefly {\em ${\rm
POSTC}_{\nu(p)}$-monoids}.\index{${\rm POSTC}_{\nu(p)}$-monoid}

We collect basic structural properties of ${\rm
POSTC}_{\nu(p)}$-monoids and show that any expanded monoid
$\mathfrak{SI}$, satisfying the following list of properties,
coincides with some ${\rm POSTC}_{\nu(p)}$-monoid
$\mathfrak{M}_{\nu(p)}$.

Let $U=U^-\,\dot{\cup}\,\{0\}\,\dot{\cup}\,U^+\,\dot{\cup}\,U'$ be
an alphabet consisting of a set $U^-$\index{$U^-$} of {\em
negative elements}\index{Element!negative}, a set
$U^+$\index{$U^+$} of {\em positive
elements}\index{Element!positive}, a set $U'$\index{$U'$} of {\em
neutral elements}\index{Element!neutral}, and zero $0$. As above,
we write $u<0$ for any element $u\in U^-$, $u>0$ for any element
$u\in U^+$, and $u\cdot v$ instead of $\{u\}\cdot\{v\}$
considering an operation $\cdot$ on the set $\mathcal{P}(U)$;
$U^{\leq 0}\rightleftharpoons U^-\cup\{0\}$\index{$U^{\leq 0}$},
$U^{\geq 0}\rightleftharpoons U^+\cup\{0\}$\index{$U^{\geq 0}$}.

A structure
$\mathfrak{M}=\langle\mathcal{P}(U);\,\cdot,\unlhd,\vee,\wedge,(\cdot\,\wedge\,\neg\,
\cdot),\circ\rangle$ is called a {\em ${\rm
POSTC}$-monoid}\index{${\rm POSTC}$-monoid} if it satisfies the
following conditions:

\medskip
${\small\bullet}$ the operation $\cdot$ of the monoid
$\langle\mathcal{P}(U)\setminus\{\varnothing\};\,\cdot\rangle$
with the unit $\{0\}$ is generated by the function $\cdot$ on
elements in $U$ such that each elements $u,v\in U$ define a
nonempty set $(u\cdot v)\subseteq U$: for any sets
$X,Y\in\mathcal{P}(U)\setminus\{\varnothing\}$ the following
equality holds:
$$
X\cdot Y=\bigcup\{u\cdot v\mid u\in X,v\in Y\};
$$
if $X\in\mathcal{P}(U)$ then $X\cdot\varnothing=\varnothing\cdot
X=\varnothing$;

\medskip
the relation $\unlhd$ on the set $\mathcal{P}(U)$ is a preorder
with the least element $\varnothing$; this preorder is induced by
the partial order $\unlhd'$ on the set $U$ of labels (forming a
upper semilattice) by the following rule: if $X,Y\in
\mathcal{P}(U)$ then $X\unlhd Y$ if and only if $X=\varnothing$,
or for any label $u\in X$ there is a label $v\in Y$ with $u\unlhd'
v$ and for any label $v\in Y$ there is a label $u\in X$ with
$u\unlhd' v$;

\medskip
${\small\bullet}$ a label $u\in U$ is called an {\em
atom}\index{Atom} if $v\unlhd u$ implies $v=u$ for any label $v\in
U$; only labels in $U^-\,\dot{\cup}\,\{0\}\,\dot{\cup}\,U^+$ may
be atoms; the label $0$ is an atom; some labels in $U^{\geq 0}$
lay under each label in $U'$, moreover, if only labels $v\in
U^{\geq 0}$ lay under a label $u\in U'$ then there is no greatest
labels among labels $v$; only labels in $U'$ lay over each label
in $U'$;

\medskip
${\small\bullet}$ the operations
$\vee,\wedge,(\cdot\,\wedge\,\neg\,\cdot)$ on the set
$U\cup\{\varnothing\}$ form a distributive lattice with relative
complements, moreover, for any elements $u,v\in
U\cup\{\varnothing\}$,
$$
u\unlhd' v\Leftrightarrow u\wedge v=u\Leftrightarrow u\vee v=v,
$$
$$
(u\wedge\neg v)=\varnothing\Leftrightarrow u\unlhd v;
$$

\medskip
${\small\bullet}$ the operation $\circ$ is defined on the set $U$
such that for any labels $u,v\in U$ the label $u\circ v$ is the
greatest element of the set $u\cdot v$;

\medskip
${\small\bullet}$ the operations $\vee,\wedge,\circ$ on the set
$\mathcal{P}(U)$ are induced by the correspondent operations on
the set $U\cup\{\varnothing\}$: if $X,Y\in \mathcal{P}(U)$ and
$\tau\in\{\vee,\wedge,\circ\}$ then $X\,\tau\, Y=\{u\,\tau\,v\mid
u\in X,v\in Y\}$; the operation $(\cdot\,\wedge\,\neg\, \cdot)$ on
the set $\mathcal{P}(U)$ is also induced by the correspondent
operation on the set $U\cup\{\varnothing\}$: if $X,Y\in
\mathcal{P}(U)$ then $X\wedge\neg Y=\{u\wedge\neg v\mid u\in
X,v\in Y\}$;

\medskip
${\small\bullet}$ the sets $U^-\cup\{\varnothing\}$ and $U^{\geq
0}\cup\{\varnothing\}$ are closed with respect to the operations
$\vee,\wedge,(\cdot\,\wedge\,\neg\,\cdot)$; the set $U'$ is closed
under the operation $\vee$; if $u\in U^-$ and $v\in U^{\geq 0}$
then $(u\vee v)\in U'$;

\medskip
${\small\bullet}$ repeating the definition in Section 2, for each
label $u\in U$, the {\em rank of semi-isolation}\index{Rank!of
semi-isolation!of label} ${\rm si}(u)\geq 1$\index{${\rm si}(u)$}
and the {\em degree of semi-isolation}\index{Degree!of
semi-isolation!of label} ${\rm deg}(u)$\index{${\rm deg}(u)$} of
label $u$ is defined inductively, ${\rm si}(\varnothing)=0$, ${\rm
deg}(\varnothing)=1$, as well as equivalence relations
$\sim_\alpha$, restrictions $X_\alpha$, $X_{\alpha,\beta}$ of sets
$X\in\{U,U\cup\{\varnothing\}\}$ and restrictions
$\mathfrak{M}'_\alpha$, $\mathfrak{M}'_{\alpha,\beta}$ for
restrictions $\mathfrak{M}'$ of $\mathfrak{M}$ to sets of labels
of ${\rm si}$-ranks $\leq\alpha$, and for labels of ${\rm
si}$-rank $\alpha$ to sets of labels of ${\rm si}$-degree
$<\beta$;

\medskip
${\small\bullet}$ the restriction
$\langle\mathcal{P}(U)\setminus\{\varnothing\};\,\cdot\rangle_{1,2}$
of the monoid
$\langle\mathcal{P}(U)\setminus\{\varnothing\};\,\cdot\rangle$ is
a $I$-groupoid;

\medskip
${\small\bullet}$ if $u<0$ then sets $u\cdot v$ and $v\cdot u$
consist of negative elements for any $v\in U$;

\medskip
${\small\bullet}$ if $u>0$ and $v>0$ then $(u\cdot v)\subseteq
U^{\geq 0}$;

\medskip
${\small\bullet}$ if $u,v\in U^{\geq 0}\cup U'$, and $u\in U'$ or
$v\in U'$, then $(u\cdot v)\subseteq U'$;

\medskip
${\small\bullet}$ for any element $u>0$ there is a nonempty set
$u^{-1}$ of {\em inverse}\index{Element!inverse} elements $u'>0$
such that $0\in (u\cdot u')\cap(u'\cdot u)$; in this case if
$u\unlhd' v$ and $v\in U^+$ then $u^{-1}\subseteq v^{-1}$;

\medskip
${\small\bullet}$ if a positive element $u$ belongs to a set
$v_1\cdot v_2$, where $v_1\circ v_2\in U^+$, then $u^{-1}\subseteq
v_2^{-1}\cdot v_1^{-1}$.

\medskip
By the definition each ${\rm POSTC}$-monoid $\mathfrak{M}$
contains ${\rm POSTC}$-submonoids $\mathfrak M^{\leq
0}$\index{$\mathfrak M^{\leq 0}$} and $\mathfrak M^{\geq
0}$\index{$\mathfrak M^{\geq 0}$} with the universes
$\mathcal{P}(U^-\cup\{0\})$ and $\mathcal{P}(U^+\cup\{0\})$
respectively, being also ${\rm POSTC}$-monoids (with $U^+\cup
U'=\varnothing$ and $U^-\cup U'=\varnothing$ respectively).

A ${\rm POSTC}$-monoid $\mathfrak{M}$ is called {\em
atomic}\index{${\rm POSTC}$-monoid!atomic} if for any label $u\in
U$ there is an atom $v\in U$ such that $v\unlhd u$.

\medskip
{\bf Theorem 6.1.} {\em For any {\rm (}at most countable and
having an ordinal ${\rm sup}\{{\rm si}(u)\mid u\in U \}${\rm )}
${\rm POSTC}$-monoid $\mathfrak{M}$ there is a {\rm (}small{\rm )}
theory $T$ with a type $p(x)\in S(T)$ and a regular labelling
function $\nu(p)$ such that $\mathfrak{M}_{\nu(p)}=\mathfrak{M}$.}

\medskip
{\em Proof} follows the same scheme as the proof of \cite[Theorem
6.1]{ShS} and, for the structure
$\langle\mathcal{P}(U)\setminus\{\varnothing\};\,\cdot\rangle_{1,2}$,
it repeats this proof word for word. Since the proof of
\cite[Theorem 6.1]{ShS} is voluminous we only point out the
distinctive features leading to the proof of this theorem.

1. A binary predicate $Q_u$ is defined for each element $u\in
U\cup\{\varnothing\}$. This predicate links only elements of the
same colors if $u\geq 0$, and defines a $Q_u$-ordered coloring
${\rm Col}$ if $u\in U^-\cup U'$; $Q_\varnothing=\varnothing$.

2.  For any elements $u,v\in U\cup\{\varnothing\}$ the following
condition is satisfied: $u\unlhd' v\Leftrightarrow Q_u\subseteq
Q_v$.

3. For any elements $u,v\in U\cup\{\varnothing\}$ the following
conditions hold:
$$u_1\vee u_2=v\Leftrightarrow Q_{u_1}\cup
Q_{u_2}=Q_v,$$
$$u_1\wedge u_2=v\Leftrightarrow Q_{u_1}\cap
Q_{u_2}=Q_v,$$
$$u_1\wedge\neg u_2=v\Leftrightarrow Q_{u_1}\setminus
Q_{u_2}=Q_v,$$
$$u_1\circ u_2=v\Leftrightarrow Q_{u_1}\circ
Q_{u_2}=Q_v.$$ In particular, the predicates $Q_{u_1}$ and
$Q_{u_2}$ are disjoint if and only if $u_1\wedge
u_2=\varnothing$.~$\Box$

\medskip
{\bf Remark 6.2.} Since labels $u\in\rho_{\nu(p)}$ for
semi-isolating formulas admit complements in $\rho_{\nu(p)}$ only
for principal types $p$ (and these complements are defined
relative to the isolating formula of $p$), unlike $I$-groupoids,
if a ${\rm POSTC}$-monoid $\mathfrak M$ is constructed by a set
$U^{\geq 0}$, it admits a representation in a transitive theory
$T$ with a (unique) type $p(x)\in S(T)$ and a regular labelling
function $\nu(p)$ such that $\mathfrak M_{\nu(p)}=\mathfrak M$ if
and only if the set-theoretic operations in $\mathfrak M$ form a
Boolean algebra.

\medskip
\centerline{\bf 7. Partial ${\rm POSTC}$-monoid} \centerline{\bf
on a set of realizations for a family of $1$-types}
\centerline{\bf of a complete theory}
\medskip

In this section, the results above for a structure of a type, as
well as results in \cite{ShS} for isolating formulas, are
generalized for a structure on a set of realizations for a family
of types.

Let $R$ be a nonempty family of types in $S^1(T)$. We denote by
$\nu(R)$\index{$\nu(R)$} a regular family of labelling functions
\index{$\rho_{\nu(R)}$}
$$\nu(p,q)\mbox{\rm : }{\rm PF}(p,q)/{\rm PE}(p,q)\to U,\,\,\,p,q\in
R,$$ $$\rho_{\nu(R)}\rightleftharpoons\bigcup\limits_{p,q\in
R}\rho_{\nu(p,q)}.$$

As in Proposition 4.1, the partial (for $|R|>1$) function ${\rm
SI}$ on the set $R\times\mathcal{P}(U)\times R$, which maps each
tuple of triples $(p_1,u_1,p_2),\ldots,(p_k,u_k,p_{k+1})$, where
$u_1\in\rho_{\nu(p_1,p_2)}\cup\{\varnothing\},\ldots,u_k\in\rho_{\nu(p_k,p_{k+1})}\cup\{\varnothing\}$,
to the set of triples $(p_1,v,p_{k+1})$, where $v\in {\rm
SI}(p_1,u_1,p_2,u_2,\ldots,p_k,u_k,p_{k+1})$, is associative:
\begin{equation}
\begin{array}{c}{\rm SI}({\rm SI}(p_1,u_1,p_2,u_2,p_3),u_3,p_4)={\rm SI}(p_1,u_1,p_2,u_2,p_3,u_3,p_4)=
\\
={\rm SI}(p_1,u_1,{\rm SI}(p_2,u_2,p_3,u_3,p_4))
\end{array}
\end{equation}
for $u_1\in\rho_{\nu(p_1,p_2)}\cup\{\varnothing\}$,
$u_2\in\rho_{\nu(p_2,p_3)}\cup\{\varnothing\}$,
$u_3\in\rho_{\nu(p_3,p_4)}\cup\{\varnothing\}$.

Consider the structure\index{$\mathfrak{M}_{\nu(R)}$}
$$\mathfrak{M}_{\nu(R)}\rightleftharpoons\langle
R\times\mathcal{P}(U)\times
R;\,\cdot,\unlhd,\vee,\wedge,(\cdot\,\wedge\,\neg\,
\cdot),\circ\rangle$$ with the partial operations $\cdot$ and
$\circ$ such that
$$(p_1,X_1,p_2)\cdot(p_2,X_2,p_3)=\bigcup\{(p_1,u_1,p_2)\cdot(p_2,u_2,p_3)\mid u_1\in X_1,u_2\in X_2\},$$
$$
(p_1,u_1,p_2)\cdot(p_2,u_2,p_3)=\{(p_1,v,p_3)\mid v\in {\rm
SI}(p_1,u_1,p_2,u_2,p_3)\},
$$
$$(p_1,X_1,p_2)\circ(p_2,X_2,p_3)=\bigcup\{(p_1,u_1,p_2)\circ(p_2,u_2,p_3)\mid u_1\in X_1,u_2\in X_2\},$$
$$
(p_1,u_1,p_2)\circ(p_2,u_2,p_3)=\{(p_1,u\circ v,p_3)\},
$$
$$u_1\in\rho_{\nu(p_1,p_2)}\cup\{\varnothing\},
u_2\in\rho_{\nu(p_2,p_3)}\cup\{\varnothing\},
$$
as well as the relation $\unlhd$ of preorder, being induced by the
partial order, of the same name, on the set of labels and the
partial operations $\vee,\wedge,(\cdot\,\wedge\,\neg\, \cdot)$
such that
$$(p,X,q)\vee(p,Y,q)=\bigcup\{(p,u,q)\vee(p,v,q)\mid u\in X,v\in Y\},$$
$$
(p,u,q)\vee(p,v,q)=\{(p,u\vee v,q)\},
$$
$$(p,X,q)\wedge(p,Y,q)=\bigcup\{(p,u,q)\wedge(p,v,q)\mid u\in X,v\in Y\},$$
$$
(p,u,q)\wedge(p,v,q)=\{(p,u\wedge v,q)\},
$$
$$(p,X,q)\wedge\neg(p,Y,q)=\bigcup\{(p,u,q)\wedge\neg(p,v,q)\mid u\in X,v\in Y\},$$
$$
(p,u,q)\wedge\neg(p,v,q)=\{(p,u\wedge\neg v,q)\},
$$
$$
u,v\in\rho_{\nu(p,q)}\cup\{\varnothing\}.
$$
The ${\rm POSTC}$-monoids $\mathfrak{M}_{\nu(p)}$, $p\in R$, are
naturally embeddable in this structure. The structure
$\mathfrak{M}_{\nu(R)}$ is called a {\em join of ${\rm
POSTC}$-monoids}\index{Join!of ${\rm POSTC}$-monoids}
$\mathfrak{M}_{\nu(p)}$, $p\in R$, relative to the family $\nu(R)$
of labelling functions and it is denoted by
$\bigoplus\limits_{\nu(R)}\mathfrak{M}_{\nu(p)}$\index{$\bigoplus\limits_{\nu(R)}\mathfrak{M}_{\nu(p)}$}.
If $\rho_{\nu(p,q)}=\varnothing$ for all $p\ne q$ the join
$\bigoplus\limits_{\nu(R)}\mathfrak{M}_{\nu(p)}$ is {\em
free},\index{Join!of ${\rm POSTC}$-monoids!free} it is represented
as the disjoint union of ${\rm POSTC}$-monoids
$\mathfrak{M}_{\nu(p)}$  and denoted by $\bigsqcup\limits_{p\in
R}\mathfrak{M}_{\nu(p)}$.\index{$\bigsqcup\limits_{p\in
R}\mathfrak{M}_{\nu(p)}$}

By (2) we have

\medskip
{\bf Proposition 7.1.} {\em For any complete theory $T$, for any
nonempty family $R\subset S(T)$ of $1$-types, and for any regular
family $\nu(R)$ of labelling functions, each $n$-ary partial
operation ${\rm
SI}(p_1,\cdot,p_2,\cdot,p_3\ldots,p_{n},\cdot,p_{n+1})$ on the set
$\mathcal{P}(U)$ is interpretable by a term of the structure
$\bigoplus\limits_{p\in\nu(R)}\mathfrak{M}_{\nu(p)}$ with fixed
types $p_1,\ldots,p_{n+1}\in R$.}

\medskip
Denote by $\mathfrak{SI}_{\nu(R)}$\index{$\mathfrak{SI}_{\nu(R)}$}
the restriction of $\mathfrak{M}_{\nu(R)}\upharpoonright
R\times(\mathcal{P}\setminus\{\varnothing\})\times R$ to the
partial operation $\cdot$.

Using Proposition 1.3 we obtain the following analogue of
Proposition 4.3.

\medskip
{\bf Proposition 7.2.} {\em For any complete theory $T$, for any
nonempty family $R\subset S(T)$ of $1$-types, and for any regular
family $\nu(R)$ of labelling functions, the restriction of the
structure $\mathfrak{SI}_{\nu(R)}$ to the set $U^{\leq 0}$ {\rm
(}respectively $U^{\geq 0}$, $U^{\geq 0}\cup U'$ {\rm )} is closed
under the partial operation $\cdot$.}

\medskip
By Proposition 7.2 the structure $\mathfrak{SI}_{\nu(R)}$ has
substructures $\mathfrak{SI}^{\leq
0}_{\nu(R)}$,\index{$\mathfrak{SI}^{\leq 0}_{\nu(R)}$}
$\mathfrak{SI}^{\geq 0}_{\nu(R)}$\index{$\mathfrak{SI}^{\geq
0}_{\nu(R)}$} and $\mathfrak{SI}^{\geq 0, {\rm
neu}}_{\nu(R)}$,\index{$\mathfrak{SI}^{\geq 0, {\rm
neu}}_{\nu(R)}$} generated by triples $(p,u,q)$ with $u\leq 0$,
$u\geq 0$, and $u\in U^{\geq 0}\cup U'$ respectively, $p,q\in R$.
Here, for any triple $(p,u,q)$ in $\mathfrak{SI}^{\geq
0}_{\nu(R)}$ the triple $(q,u^{-1},p)$ is also attributed to
$\mathfrak{SI}^{\geq 0}_{\nu(R)}$.

\medskip
Replacing for the definition in Section 5 the function $\nu(p)$ to
the family $\nu(R)$ of functions we obtain the notions of {\em
$(\alpha,n)$-deterministic, $\alpha$-deterministic, almost
$\alpha$-deterministic}, and {\em $(\alpha,\omega)$-deterministic}
structures $\mathfrak{SI}_{\nu(p)}\upharpoonright
U_0$.\index{Structure!$(\alpha,n)$-deterministic}\index{Structure!$\alpha$-deterministic}\index{Structure!almost
$\alpha$-deterministic}\index{Structure!$(\alpha,\omega)$-deterministic}

Below we formulate a series of assertions that immediately
transform from the class of structures $\mathfrak{SI}_{\nu(p)}$ to
the class of structures $\mathfrak{SI}_{\nu(R)}$.

\medskip
{\bf Proposition 7.3.} {\em {\rm (Monotony)\index{Monotony}} If a
structure $\mathfrak{SI}_{\nu(R)}\upharpoonright U_0$ is {\rm
(}almost{\rm )} $\alpha$-deterministic and $\beta$ is a positive
ordinal then the structure $(\mathfrak{SI}_{\nu(R)}\upharpoonright
U_0)\upharpoonright\beta$ is also {\rm (}almost{\rm )}
$\alpha$-deterministic.}

\medskip
{\bf Proposition 7.4.} {\em For a structure
$\mathfrak{SI}_{\nu(R)}$ and ordinals $\alpha,\beta$, where
$\alpha,\beta>0$, $\beta\in\omega+1$, the following conditions are
equivalent:

{\rm (1)} the structure $\mathfrak{SI}_{\nu(R)}$ is
$(\alpha,\beta)$-deterministic;

{\rm (2)} for any types $p,q,r\in R$ and labels
$u_1\in\rho_{\nu(p,q),\alpha,\beta}$,
$u_2\in\rho_{\nu(q,r),\alpha,\beta}$ the inequality ${\rm
si}(u_1\circ u_2)\leq\alpha$ holds and is ${\rm si}(u_1\circ
u_2)=\alpha$ then ${\rm deg}(u_1\circ u_2)<\beta$.}

\medskip
{\bf Corollary 7.5.} {\em For a structure $\mathfrak{SI}_{\nu(R)}$
and a positive ordinal $\alpha$, the following conditions are
equivalent:

{\rm (1)} the structure $\mathfrak{SI}_{\nu(R)}$ is
$\alpha$-deterministic;

{\rm (2)} ${\rm si}(u_1\circ u_2)\leq\alpha$ for any types
$p,q,r\in R$ and labels $u_1\in\rho_{\nu(p,q),\alpha}$,
$u_2\in\rho_{\nu(q,r),\alpha}$ .}

\medskip
{\bf Corollary 7.6.} {\em If ${\rm si}(R)$ is an ordinal then the
structure $\mathfrak{SI}_{\nu(R)}$ is almost ${\rm
si}(R)$-deterministic.}

\medskip
{\bf Proposition 7.7.} {\em If a structure
$\mathfrak{SI}_{\nu(R)}$ is $(\alpha,\beta)$-deterministic then
$\mathfrak{SI}_{\nu(R),\alpha}\rightleftharpoons
\mathfrak{SI}_{\nu(R)}\upharpoonright\alpha$\index{$\mathfrak{SI}_{\nu(R),\alpha}$}
is also $(\alpha,\beta)$-deterministic.}

\medskip
{\bf Proposition 7.8.} {\em If ${\rm si}(R)$ is an ordinal then
the structure $\mathfrak{SI}_{\nu(R)}$ is ${\rm
si}(R)$-deterministic if and only if the value ${\rm deg}(R)$ is
not defined or equals $1$.}

\medskip
{\bf Proposition 7.9.} {\em A structure $\mathfrak{P}_{\nu(R)}$ is
{\rm (}almost{\rm )} deterministic if and only if
$\mathfrak{SI}_{\nu(R),1,2}$ is {\rm (}almost{\rm )}
$1$-deterministic.}

\medskip
Let $R$ be a nonempty family of complete $1$-types of a theory
$T$, $\nu(R)$ be a regular family of labelling functions, and
$\alpha$ be an ordinal, $\alpha>0$. The structure
$\mathfrak{SI}_{\nu(R)}$ is called {\em locally
$\alpha$-deterministic}\index{Structure!locally
$\alpha$-deterministic} if for any nonempty finite set
$R_0\subseteq R$ there is a natural number $n\geq 2$ such that the
structure $\mathfrak{SI}_{\nu(R_0)}$ is
$(\alpha,n)$-deterministic.

Repeating the proof of Proposition 5.8 we obtain

\medskip
{\bf Proposition 7.10.} {\em Let $R$ be a nonempty family of
complete $1$-types of a theory $T$, $\nu(R)$ be a regular family
of labelling functions, ${\rm si}(R)<\omega$. The following
conditions are equivalent:

$(1)$ the structure $\mathfrak{SI}_{\nu(R)}$ is locally
$1$-deterministic;

$(2)$ the set $\rho_{\nu(p,q)}$ is finite for any $p,q\in R$;

$(3)$ the set $\rho_{\nu(p,q),1}$ is finite for any $p,q\in R$;

$(4)$ the set $\rho_{\nu(p,q),1,2}$ {\rm (}consisting of all atoms
$u\in\rho_{\nu(p,q)}${\rm )} is finite for any $p,q\in R$.}

\medskip
The notion of $(n,p)$-type is generalized in the following
definition.

\medskip
{\bf Definition} (K.~Ikeda, A.~Pillay, A.~Tsuboi \cite{IPT}). Let
$p_1(x_1),\ldots,$ $p_n(x_n)$ be types in $S(T)$ with disjoint
free variables. A type $q(x_1,\ldots,x_n)\in S(T)$ is said to be a
\emph{$(p_1,\ldots,p_n)$-type}\index{$(p_1,\ldots,p_n)$-type} if
$q(x_1,\ldots,x_n)\supseteq\bigcup\limits_{i=1}^n p_i(x_i)$. The
set  of  all $(p_1,\ldots,p_n)$-types  of $T$  is denoted
by~$S_{p_1,\ldots,p_n}(T)$\index{$S_{p_1,\ldots,p_n}(T)$}. A
theory $T$ is {\em almost
$\omega$-categorical}\index{Theory!almost $\omega$-categorical} if
for any types $p_1(x_1),\ldots,p_n(x_n)\in S(T)$ there are only
finitely many types $q(x_1,\ldots,x_n)\in S_{p_1,\ldots,p_n}(T)$.

\medskip
{\bf Definition} (B.~S.~Baizhanov\index{Baizhanov B. S.},
S.~V.~Sudoplatov\index{Sudoplatov S. V.},
V.~V.~Verbovskiy\index{Verbovskiy V. V.} \cite{BSV}). A type
$q(\bar{x})$ in $S_{p_1,\ldots,p_n}(T)$ is said to be
\emph{$(p_1,\ldots,p_n)$-principal}\index{Type!$(p_1,\ldots,p_n)$-principal}
if there is a formula $\varphi(\bar{y})\in q(\bar{x})$ such that
$$\cup\{p_i(x_i)\mid
i=1,\ldots,n\}\cup\{\varphi(\bar{x})\}\vdash q(\bar{x}).$$

\medskip
The following lemma obviously generalizes Lemma 5.10.

\medskip
{\bf Lemma 7.11} \cite{BSV}. {\em For any types
$p_1(x_1),\ldots,p_n(x_n)\in S(\varnothing)$ the following
conditions are equivalent:

{\rm (1)} the set of $(p_1,\ldots,p_n)$-types with free variables
in $(x_1,\ldots,x_n)$ is finite;

{\rm (2)} any $(p_1,\ldots,p_n)$-type is
$(p_1,\ldots,p_n)$-principal.}

\medskip
By Lemma 7.11, a theory $T$ is almost $\omega$-categorical if and
only if for any types $p_1(x_1),\ldots,p_n(x_n)\in S^1(T)$, each
$(p_1,\ldots,p_n)$-type is $(p_1,\ldots,p_n)$-principal.

\medskip
Proposition 7.10 and Lemma 7.11 imply

\medskip
{\bf Corollary 7.12.} {\em If $R$ is a nonempty family of complete
$1$-types of a theory $T$, $\nu(R)$ is a regular family of
labelling functions, and all $(p_1,p_2)$-types, where $p_1,p_2\in
R$, are $(p_1,p_2)$-principal then the structure
$\mathfrak{SI}_{\nu(R)}$ is locally $1$-deterministic.}

\medskip
{\bf Corollary 7.13.} {\em If $T$ is an almost
$\omega$-categorical theory and $\nu(S^1(\varnothing))$ is a
regular family of labelling functions then the structure
$\mathfrak{SI}_{\nu(R)}$ is locally $1$-deterministic.}

\medskip
For a nonempty family $R$ of $1$-types in $S(T)$ and a positive
ordinal $\alpha$, we denote by ${\rm SI}_{R,\alpha}$\index{${\rm
SI}_{R,\alpha}$} (in a model $\mathcal{M}$ of $T$) the restriction
of ${\rm SI}_R$ to the set of formulas of ${\rm si}$-ranks
$\leq\alpha$:
$$
{\rm SI}_{R,\alpha}\rightleftharpoons\{(a,b)\mid {\rm tp}(a),{\rm
tp}(b)\in R\mbox{ and }a\mbox{ semi-isolates }b$$
$$\mbox{ by a formula }\theta_{{\rm tp}(a),u,{\rm
tp}(b)}(x,y),\mbox{ with a }{\rm si}\mbox{-rank }\leq\alpha\}.
$$

Clearly, $I_R={\rm SI}_{R,1}$ for any nonempty family $R$ of
$1$-types. Considering this equality and the equality
$\mathfrak{SI}_{\nu(R),1,2}=\mathfrak{P}_{\nu(R)}$, the following
proposition generalizes Proposition 5.13 as well as Propositions
4.3 and 8.3 in \cite{ShS}.

\medskip
{\bf Proposition 7.14.} {\em Let $R$ be a nonempty family of
complete $1$-types of a theory $T$, $\nu(R)$ be a regular family
of labelling functions, and $\alpha$ be a positive ordinal. The
following conditions are equivalent:

$(1)$ the relation ${\rm SI}_{R,\alpha}$ {\rm (}on a set of
realizations of types $p\in R$ in any model $\mathcal{M}\models
T${\rm )} is transitive;

$(2)$ the structure $\mathfrak{SI}_{\nu(R),\alpha}$ is almost
$\alpha$-deterministic.}

\medskip
{\em Proof} repeats the proof of Proposition 5.13 almost word for
word.~$\Box$

\medskip
Propositions 1.3 and 7.14 imply the following assertions.

\medskip
{\bf Corollary 7.15.} {\em Let $R$ be a nonempty family of
complete $1$-types of a theory $T$, $\nu(R)$ be a regular family
of labelling functions, and $\alpha$ be a positive ordinal. The
following conditions are equivalent:

$(1)$ the relation ${\rm SI}_{R,\alpha}$ in any model
$\mathcal{M}\models T$ is a partial order;

$(2)$ the structure $\mathfrak{SI}_{\nu(R),\alpha}$ is almost
$\alpha$-deterministic and $\rho_{\nu(R),\alpha}\subseteq U^{\leq
0}$.

The partial order ${\rm SI}_{R,\alpha}$ is identical if and only
if $\rho_{\nu(R),\alpha}=\{0\}$. The non-identical partial order
${\rm SI}_{R,\alpha}$ has infinite chains if and only if
$|\rho_{\nu(p),\alpha}|>1$ for some type $p\in R$ or there is a
sequence $p_n$, $n\in\omega$, of pairwise distinct types in $R$
such that $|\rho_{\nu(p_n,p_{n+1}),\alpha}|\geq 1$, $n\in\omega$,
or $|\rho_{\nu(p_{n+1},p_{n})}|\geq 1$, $n\in\omega$.}

\medskip
{\bf Corollary 7.16.} {\em Let $R$ be a nonempty family of
complete $1$-types of a theory $T$, $\nu(R)$ be a regular family
of labelling functions, and $\alpha$ be a positive ordinal. The
following conditions are equivalent:

$(1)$ the relation ${\rm SI}_{R,\alpha}$ on a set of realizations
of types $p\in R$ in any model $\mathcal{M}\models T$ is an
equivalence relation;

$(2)$ the structure $\mathfrak{SI}_{\nu(R),\alpha}$ is almost
$\alpha$-deterministic and $\rho_{\nu(R),\alpha}\subseteq U^{\geq
0}$.}

\medskip
The results above substantiate that the diagram in Figure 1 admits
the transformation replacing the type $p$ by a nonempty family
$R\subseteq S^1(\varnothing)$.

\medskip
\centerline{\bf 8. ${\rm POSTC}_\mathcal{R}$-structures}
\medskip

{\bf Definition.} Let $\mathcal{R}$ be a nonempty set,
$$U=U^-\,\dot{\cup}\,\{0\}\,\dot{\cup}\,U^+\,\dot{\cup}\,U'$$
be an alphabet consisting of a set $U^-$\index{$U^-$} of {\em
negative elements}\index{Element!negative}, a set
$U^+$\index{$U^+$} of {\em positive
elements}\index{Element!positive}, a set $U'$\index{$U'$} of {\em
neutral elements}\index{Element!neutral}, and zero $0$. If $p$ and
$q$ are elements in $\mathcal{R}$, we write $u<0$ and $(p,u,q)<0$
for any element $u\in U^-$, $u>0$ and $(p,u,q)>0$ for any element
$u\in U^+$; $U^{\leq 0}\rightleftharpoons
U^-\cup\{0\}$\index{$U^{\leq 0}$}, $U^{\geq 0}\rightleftharpoons
U^+\cup\{0\}$\index{$U^{\geq 0}$}. For the set $\mathcal{R}^2$ of
all pairs $(p,q)$, $p,q\in \mathcal{R}$, we consider a {\em
regular} family\index{Family!of sets!regular}
$\mu(\mathcal{R})$\index{$\mu(\mathcal{R})$} of sets
$\mu(p,q)\subseteq U$ such that

\medskip
${\small\bullet}$ $0\in\mu(p,q)$ if and only if $p=q$;

\medskip
${\small\bullet}$ $\mu(p,p)\cap\mu(q,q)=\{0\}$ for $p\ne q$;

\medskip
${\small\bullet}$ $\mu(p,q)\cap\mu(p',q')=\varnothing$ if $p\ne q$
and $(p,q)\ne(p',q')$;

\medskip
${\small\bullet}$ $\bigcup\limits_{p,q\in\mathcal{R}}\mu(p,q)=U$.

\medskip
Further we write $\mu(p)$ instead of $\mu(p,p)$, and considering a
partial operation $\cdot$ on the set
$\mathcal{R}\times\mathcal{P}(U)\times\mathcal{R}$ we shall write,
as above, $(p,u,q)\cdot(q,v,r)$ instead of
$(p,\{u\},q)\cdot(q,\{v\},r)$.

A structure
$$\mathfrak{M}=\langle\mathcal{R}\times\mathcal{P}(U)\times\mathcal{R};\,\cdot,\unlhd,\vee,\wedge,(\cdot\,\wedge\,\neg\,
\cdot),\circ\rangle$$ with a regular family $\mu(\mathcal{R})$ of
sets is said to be a {\em ${\rm
POSTC}_\mathcal{R}$-structure}\index{${\rm
POSTC}_\mathcal{R}$-structure} if the following conditions hold:

\medskip
${\small\bullet}$ the partial operation $\cdot$ of the structure
$\langle
\mathcal{R}\times\mathcal{P}(U)\times\mathcal{R};\,\cdot\rangle$
has values $(p,X,q)\cdot(p',Y,q')$ only for $p'=q$,
$X\subseteq\mu(p,q)$, $Y\subseteq\mu(p',q')$, and it is generated
by the function $\cdot$ for elements in $U$: for any sets
$X,Y\in\mathcal{P}(U)$, $\varnothing\ne X\subseteq\mu(p,q)$,
$\varnothing\ne Y\subseteq\mu(q,r)$, the following equality is
satisfied:
$$
(p,X,q)\cdot(q,Y,r)=\bigcup\{(p,x,q)\cdot(q,y,r)\mid x\in X,y\in
Y\},
$$
and if some of $X,Y$ is empty then
$(p,X,q)\cdot(q,Y,r)=\varnothing$;

\medskip
${\small\bullet}$ each restriction $\mathfrak
M_{\mu(p)}$\index{$\mathfrak M_{\mu(p)}$} of $\mathfrak M$ to
$\{p\}\times\mathcal{P}(\mu(p))\times\{p\}$ is isomorphic to a
${\rm POSTC}$-monoid with the universe $\mathcal{P}(\mu(p))$,
$p\in\mathcal{R}$; atoms $u\in\mu(p)$ in $\mathfrak M_{\mu(p)}$
are called {\em $p$-atoms};\index{$p$-atom}

\medskip
${\small\bullet}$ each restriction $\mathfrak
M_{\mu(p,q)}$\index{$\mathfrak M_{\mu(p,q)}$}, $p\ne q$, of
$\mathfrak M$ to $\{p\}\times\mathcal{P}(\mu(p,q))\times\{q\}$ has
empty partial operations $\cdot$ and $\circ$; the restriction of
$\mathfrak M_{\mu(p,q)}$ to the relation $\unlhd$ is a preordered
set $\langle
\{p\}\times\mathcal{P}(\mu(p,q))\times\{q\};\,\unlhd_{p,q}\rangle$
with the least element $(p,\varnothing,q)$, the preorder
$\unlhd_{p,q}$ of this structure is induced by the partial order
$\unlhd'_{p,q}$ on the set $\mu(p,q)$ of labels (forming a upper
semilattice if $\mu(p,q)\ne\varnothing$) by the following rule: if
$X,Y\in \mathcal{P}(\mu(p,q))$ then $X\unlhd_{p,q} Y$ if and only
if $X=\varnothing$, or for any label $u\in X$ there is a label
$v\in Y$ with $u\unlhd_{p,q} v$ and for any label $v\in Y$ there
is a label $u\in X$ with $u\unlhd_{p,q} v$;

\medskip
${\small\bullet}$ a label $u\in \mu(p,q)$, where $p\ne q$, is said
to be a {\em $(p,q)$-atom}\index{$(p,q)$-atom} if $v\unlhd_{p,q}u$
implies $v=u$ for any label $v\in \mu(p,q)$; only labels in
$\mu(p,q)\cap(U^-\cup U^+)$ may be $(p,q)$-atoms; some labels in
$\mu(p,q)\cap U^{\geq 0}$ lay under each label in $\mu(p,q)\cap
U'$, moreover, if only labels $v\in\mu(p,q)\cap U^{\geq 0}$ lay
under a label $u\in \mu(p,q)\cap U'$ then there are no greatest
labels among labels $v$; only labels in $\mu(p,q)\cap U'$ lay over
each label in $\mu(p,q)\cap U'$;

\medskip
${\small\bullet}$ the operations
$\vee,\wedge,(\cdot\,\wedge\,\neg\,\cdot)$ are defined on each set
$\mu(p,q)\cup\{\varnothing\}$ in the structure $\mathfrak
M_{\mu(p,q)}$ and form a distributive lattice with relative
complements on $\mu(p,q)\cup\{\varnothing\}$, moreover, for any
elements $u,v\in \mu(p,q)\cup\{\varnothing\}$,
$$
u\unlhd_{p,q} v\Leftrightarrow u\wedge v=u\Leftrightarrow u\vee
v=v\Leftrightarrow u\wedge\neg v=\varnothing;
$$

\medskip
${\small\bullet}$ the relation $\unlhd$ on the set
$\mathcal{R}\times\mathcal{P}(U)\times\mathcal{R}$ is a preorder
with minimal elements $(p,\varnothing,q)$, $p,q\in\mathcal{R}$;
this preorder is induced by the union $\unlhd_U$ of preorders
$\unlhd_p$ in the structures $\mathfrak M_{\mu(p)}$,
$p\in\mathcal{R}$, and of preorders $\unlhd_{p,q}$ in the
structures $\mathfrak M_{\mu(p,q)}$, $p,q\in\mathcal{R}$, $p\ne
q$, on sets of labels in these structures: if $X,Y\in
\mathcal{P}(U)$ then $(p,X,q)\unlhd(p',Y,q')$ if and only if
$p=p'$, $q=q'$, and $X=\varnothing$ or for any label $u\in X$
there is a label $v\in Y$ with $u\unlhd_U v$ and for any label
$v\in Y$ there is a label $u\in X$ with $u\unlhd_U v$;

\medskip
${\small\bullet}$ the partial operations
$\vee,\wedge,(\cdot\,\wedge\,\neg\,\cdot)$ are defined on the set
$\mathcal{R}\times (U\cup\{\varnothing\})\times\mathcal{R}$ in the
structure $\mathfrak{M}$ being unions of correspondent operations
on the sets $\mu(p)\cup\{\varnothing\}$ in $\mathfrak M_{\mu(p)}$
and on the sets $\mu(p,q)\cup\{\varnothing\}$ in $\mathfrak
M_{\mu(p,q)}$, $p\ne q$;

\medskip
${\small\bullet}$ the partial operation $\circ$ is defined on the
set $\mathcal{R}\times (U\cup\{\varnothing\})\times\mathcal{R}$ in
the structure $\mathfrak{M}$ being obtained from the union of
correspondent operations in the structures $\mathfrak M_{\mu(p)}$,
$p\in\mathcal{R}$, by the following extension: if $u_1\in\mu(p,q)$
and $u_2\in\mu(q,r)$ then there is unique element $v\in\mu(p,r)$,
such that $(p,u_1,q)\circ(q,u_2,r)=(p,v,r)$; this element $v$ is
the $\unlhd_{p,r}$-greatest label in the set
$(p,u_1,q)\cdot(q,u_2,r)$, it is called a {\em
composition}\index{Composition!of elements} of elements $u_1$ and
$u_2$ and it is denoted by $u_1\circ u_2$;
$$(p,u_1,q)\circ(q,\varnothing,r)=(p,\varnothing,q)\cdot(q,u_2,r)=(p,\varnothing,q)\cdot(q,\varnothing,r)=(p,\varnothing,r);$$

\medskip
${\small\bullet}$ the partial operations $\vee,\wedge,\circ$ on
the set $\mathcal{R}\times\mathcal{P}(U)\times\mathcal{R}$ are
induced by the correspondent partial operations on the set
$\mathcal{R}\times(U\cup\{\varnothing\})\times\mathcal{R}$: if
$(p,X,q),(p',Y,q')\in
\mathcal{R}\times\mathcal{P}(U)\times\mathcal{R}$ and
$\tau\in\{\vee,\wedge,\circ\}$ then the value
$(p,X,q)\,\tau\,(p',Y,q')$ is not defined or it is defined and
coincides with the set $\{(p,u,q)\,\tau\,(p',v,q')\mid u\in X,v\in
Y\}$, in which all values are defined; the partial operation
$(\cdot\,\wedge\,\neg\, \cdot)$ on the set
$\mathcal{R}\times\mathcal{P}(U)\times\mathcal{R}$ is also induced
by the correspondent partial operation on the set
$\mathcal{R}\times(U\cup\{\varnothing\})\times\mathcal{R}$: if
$(p,X,q),(p',Y,q')\in
\mathcal{R}\times\mathcal{P}(U)\times\mathcal{R}$ then the value
$(p,X,q)\wedge\neg (p',Y,q')$ is defined only for $p=p'$, $q=q'$,
$X,Y\subseteq\mu(p,q)$ and it is equal to
$\{(p,u,q)\wedge\neg(p,v,q)\mid u\in X,v\in Y\}$;

\medskip
${\small\bullet}$ each of the sets $U^-\cup\{\varnothing\}$ and
$U^{\geq 0}\cup\{\varnothing\}$ is closed under operations
$\vee,\wedge,(\cdot\,\wedge\,\neg\,\cdot)$; the set $U'$ is closed
under the operation $\vee$; if $u\in U^-$ and $v\in U^{\geq 0}$
then $(u\vee v)\in U'$;

\medskip
${\small\bullet}$ repeating the definition in Section 2, each
label $u\in U$ obtains inductively the {\em rank of
semi-isolation}\index{Rank!of semi-isolation!of label} ${\rm
si}(u)\geq 1$\index{${\rm si}(u)$} and the {\em degree of
semi-isolation}\index{Degree!of semi-isolation!of label} ${\rm
deg}(u)$\index{${\rm deg}(u)$}, ${\rm si}(\varnothing)=0$, ${\rm
deg}(\varnothing)=1$, as well as the following attributes are
defined: the equivalence relations $\sim_\alpha$, restrictions
$X_\alpha$ and $X_{\alpha,\beta}$ of sets
$X\in\{U,U\cup\{\varnothing\}\}$, and restrictions
$\mathfrak{M}'_\alpha$, $\mathfrak{M}'_{\alpha,\beta}$ for
restrictions $\mathfrak{M}'$ of the structure $\mathfrak{M}$ to
the set of labels of ${\rm si}$-ranks $\leq\alpha$, and for labels
of ${\rm si}$-rank $\alpha$ to the set of labels of ${\rm
si}$-degree $<\beta$;

\medskip
${\small\bullet}$ the restriction
$\langle\mathcal{R}\times(\mathcal{P}(U)\setminus\{\varnothing\})\times\mathcal{R};\,\cdot\rangle_{1,2}$
of the structure $\mathfrak{M}$ is an $I_\mathcal{R}$-structure;

\medskip
${\small\bullet}$ if $u\in\mu(p,q)$ and $u<0$ then the set
$(p,u,q\cdot(q,v,r)$ and $(r,v',p)\cdot(p,u,q)$ consist of
negative elements for any $v\in \mu(q,r)$ and $v'\in(r,p)$;

\medskip
${\small\bullet}$ if $u\in\mu(p,q)$, $v\in\mu(q,r)$, $u>0$, and
$v>0$, then the set $(p,u,q)\cdot(q,v,r)$ consists of elements in
$U^{\geq 0}$;

\medskip
${\small\bullet}$ if $u\in\mu(p,q)\cap(U^{\geq 0}\cup U')$,
$v\in\mu(q,r)\cap(U^{\geq 0}\cup U')$, and $u\in U'$ or $v\in U'$,
then $(p,u,q)\cdot(q,v,r)\subseteq U'$;

\medskip
${\small\bullet}$ for any element $u\in\mu(p,q)$ with $u>0$ there
is a nonempty set $u^{-1}$ of {\em inverse}\index{Element!inverse}
elements $u'>0$ such that $(p,0,p)\in (p,u,q)\cdot(q,u',p)$ and
$(q,0,q)\in(q,u',p)\cdot(p,u,q)$, moreover, if $u\unlhd'_{p,q} v$
and $v\in U^+$ then $u^{-1}\subseteq v^{-1}$;

\medskip
${\small\bullet}$ if an element $(p,u,r)$, where $u>0$, belongs to
a set $(p,v_1,q)\cdot(q,v_2,r)$, where $v_1\circ v_2\in U^+$, then
$(r,u^{-1},p)\subseteq(r,v_2^{-1},q)\cdot(q,v_1^{-1},p)$.

\medskip
By the definition, each ${\rm POSTC}_\mathcal{R}$-structure
$\mathfrak{M}$ contains ${\rm POSTC}_\mathcal{R}$-substruc\-tures
$\mathfrak M^{\leq 0}$\index{$\mathfrak M^{\leq 0}$} and
$\mathfrak M^{\geq 0}$\index{$\mathfrak M^{\geq 0}$} being
restrictions of $\mathfrak M$ to the sets $U^{\leq 0}$ and
$U^{\geq 0}$ respectively.

A ${\rm POSTC}_\mathcal{R}$-structure $\mathfrak{M}$ is called
{\em atomic}\index{${\rm POSTC}$-structure!atomic} if for any
label $u\in\mu(p)$, $p\in\mathcal{R}$, there is a $p$-atom $v\in
U$ such that $v\unlhd_p u$, and for any label $u\in\mu(p,q)$,
$p,q\in\mathcal{R}$, $p\ne q$, there is a $(p,q)$-atom $v\in U$
such that $v\unlhd_{p,q} u$.

Combining the proof of Theorems 6.1 and 9.1 in \cite{ShS} as well
as the proof of Theorem 6.1, we obtain the following theorem.

\medskip
{\bf Theorem 8.1.} {\em For any \ {\rm (}at most countable and
having an ordinal \ ${\rm sup}\{{\rm si}(u)\mid u\in U \}${\rm
)}{\rm )} ${\rm POSTC}_\mathcal{R}$-structure $\mathfrak{M}$ there
is a {\rm (}small{\rm )} theory $T$ with a family of $1$-types
$R\subset S(T)$ and a regular family $\nu(R)$ of labelling
functions such that $\mathfrak{M}_{\nu(R)}=\mathfrak{M}$.}

\medskip
In conclusion, we note that, using the operation $\cdot^{\rm eq}$,
the constructions above can be transformed for an arbitrary family
of types in $S(T)$.

\end{document}